\documentclass[final,hidelinks,onefignum,onetabnum]{siamart220329}

\usepackage{tikz}
\usetikzlibrary{shapes}
\tikzset{My Arrow Style/.style={single arrow, fill=red!50, anchor=base, align=center,text width=2.8cm}}

\usepackage{esint}
\usepackage{xcolor}
\usepackage[export]{adjustbox}

\makeatletter

\usepackage{amsfonts}
\usepackage{subfigure}
\usepackage{epsfig}
\usepackage[small,compact]{titlesec}
\usepackage{paralist}
 {\usepackage{color}}{}
\usepackage{caption}
\usepackage{pdflscape}
\usepackage{multirow}
\usepackage{enumitem}
\usepackage{lscape}
\usepackage{lineno}

\newtheorem{principle}{Principle}

\usepackage{lipsum}
\usepackage{amsfonts}
\usepackage{graphicx}
\usepackage{epstopdf}
\usepackage{algorithmic}
\ifpdf
  \DeclareGraphicsExtensions{.eps,.pdf,.png,.jpg}
\else
  \DeclareGraphicsExtensions{.eps}
\fi

\newsiamremark{remark}{Remark}
\newsiamremark{hypothesis}{Hypothesis}
\crefname{hypothesis}{Hypothesis}{Hypotheses}
\newsiamthm{claim}{Claim}

\headers{A meshfree RBF-FD constant along normal method for solving PDEs on surfaces}{V\'ictor Bayona, Argyrios Petras, C\'ecile Piret, and Steven J. Ruuth}

\title{A meshfree RBF-FD constant along normal method for solving PDEs on surfaces}

\author{V\'ictor Bayona\thanks{Departamento de Matem\'aticas, Universidad Carlos III de Madrid, Legan\'es 28911, Spain 
  (\email{vbayona@math.uc3m.es}).}
\and Argyrios Petras\thanks{RICAM - Johann Radon Institute for Computational and Applied Mathematics, Austrian Academy of Sciences, Linz 4040, Austria 
  (\email{argyrios.petras@ricam.oeaw.ac.at}.)}
\and C\'ecile Piret\thanks{Department of Mathematical Sciences, Michigan Technological University, Houghton, MI 49931, USA
  (\email{cmpiret@mtu.edu}.)}
\and Steven J. Ruuth\thanks{Department of Mathematics, Simon Fraser University, Burnaby, V5A 1S6, BC, Canada 
  (\email{sruuth@sfu.ca}).}
}
\usepackage{amsopn}

\ifpdf
\hypersetup{
  pdftitle={A meshfree RBF--FD constant along normal method for solving PDEs on surfaces},
  pdfauthor={V\'ictor Bayona, Argyrios Petras, C\'ecile Piret, and Steven J. Ruuth}
}
\fi

\begin{document}

\nolinenumbers

\maketitle

\begin{abstract}
This paper introduces a novel meshfree methodology based on Radial Basis Function-Finite Difference (RBF-FD) approximations for the numerical solution of partial differential equations (PDEs) on surfaces of codimension 1 embedded in $\mathbb{R}^3$. The method is built upon the principles of the closest point method, without the use of a grid or a closest point mapping. We show that the combination of local embedded stencils with these principles can be employed to approximate surface derivatives using polyharmonic spline kernels and polynomials (PHS+Poly) RBF-FD. Specifically, we show that it is enough to consider a constant extension along the normal direction only at a single node to overcome the rank deficiency of the polynomial basis. An extensive parameter analysis is presented to test the dependence of the approach. We demonstrate high-order convergence rates on problems involving surface advection and surface diffusion, and solve Turing pattern formations on surfaces defined either implicitly or by point clouds. Moreover, a simple coupling approach with a particle tracking method demonstrates the potential of the proposed method in solving PDEs on evolving surfaces in the normal direction. Our numerical results confirm the stability, flexibility, and high-order algebraic convergence of the approach.
\end{abstract}



\section{Introduction}
PDEs on surfaces appear throughout the natural and applied sciences. There are many methods for the numerical solution of such PDEs, among which meshfree methods \cite{chen2020extrinsic,fuselier2013high,jones2023generalized,lehto2017radial,piret2012orthogonal,shankar2018rbf,shankar2015radial,
Suchde2019} and the closest point method \cite{macdonald2009implicit,petras2018rbf,ruuth2008simple} have achieved a significant attention over the last years. Although these two approaches are conceptually very different, both formulate differential operators on surfaces entirely in Cartesian (or extrinsic) coordinates, avoiding any singularities associated with the use of intrinsic coordinate systems.


One common approach in meshfree methods is known as the the \textit{projected gradient method}, first introduced in \cite{flyer2009radial} using a global RBF interpolation on the sphere and extended later to radial basis function generated finite differences (RBF-FD) in \cite{flyer2012guide}. This was later generalized to arbitrary surfaces using global RBFs in \cite{fuselier2013high}, and using RBF-FD in \cite{lehto2017radial,shankar2015radial}. The underlying idea in this approach is to use RBF interpolants to approximate first order differential operators such as the surface gradient 
or the surface divergence, and then represents higher order operators as compositions of the approximations of the first order operators, such as the Laplace--Beltrami operator, i.e. $\Delta_{\Gamma}=\nabla_{\Gamma}\,\cdot\,\nabla_{\Gamma}$. 

In all these works, infinitely smooth RBFs, such as the inverse multiquadric or Gaussians, are used, which depend on a shape parameter. For smooth target functions, smaller values of this parameter generally lead to more accurate RBF interpolants. However, the standard numerical approach becomes ill-conditioned for small shape parameter values. As a result, its selection 
limits the efficiency and accuracy of the approach. Stabilization techniques have been proposed to tackle the ill-conditioning for Gaussian RBFs \cite{fornberg_flyer_2015}; however, they can be computationally expensive.

To overcome this issue, the works in \cite{BFFB17,FFBB16} proposed an alternative RBF-FD formulation based on poly-harmonic splines (PHS) RBFs augmented with polynomials (PHS+poly). In this formulation, the RBFs are independent of a shape parameter, and the convergence rate is determined by the degree of the augmented polynomials. Unfortunately, the application of this approach to the numerical solution of surface differential equations is not straightforward. The main problem is that the addition of polynomial basis terms to the RBF interpolant leads to rank deficient or ill-conditioned interpolation matrices when interpolating functions on manifolds due to linear dependence of polynomials.

Different approaches have been used to address this issue. In \cite{shankar2018rbf} the polynomial basis elements augmenting the PHS are obtained using the least orthogonal interpolation technique (LOI) on each RBF-FD stencil to obtain local restrictions of polynomials in $\mathbb{R}^3$ to stencils on a surface. This approach was extended to advection problems on arbitrary surfaces in \cite{shankar2020robust}.   

An alternative method based on a different philosophy was proposed in \cite{reeger2016numerical,reeger2018numerical} for numerical quadrature over smooth surfaces. The idea is to project the RBF-FD stencil nodes to the tangent plane, compute the RBF-FD weights there using PHS plus 2D polynomials, and then transform the weights back to the original stencil location. This approach was extended in \cite{gunderman2020transport} to solve transport equations over the sphere using Householder reflections. 

One of the procedures developed outside the RBF community follows independently the same idea, known as the \textit{tangent plane method}. This approach was first introduced in \cite{demanet2006painless} for approximating the surface Laplacian using polynomial based approximations, and generalized later to other surface differential equations using polynomial weighted least squares in \cite{Suchde2019}. This approach was first considered in \cite{shaw2019radial} for approximating the Laplace-Beltrami operator with RBF-FD, and very recently was presented in \cite{jones2023generalized} a RBF-FD method for approximating the tangent space of surfaces defined only by point clouds.

The philosophy behind the closest point method is somehow different. This requires a closest point representation of a surface and the use of two main principles to decouple the geometry from the PDE formulation: the \emph{equivalence of gradients} and the \emph{equivalence of divergence}.
\begin{principle}\label{p1}
Let $u$ be a function on $\mathbb{R}^d$ that is constant-along-normal directions of a surface $\Gamma$. Then, intrinsic gradients are equivalent to standard gradients at the surface, $\nabla_\Gamma u = \nabla u$.
\end{principle}
\begin{principle}\label{p2}
Let $\boldsymbol{v}$ be a vector field on $\mathbb{R}^d$ that is tangent to a surface $\Gamma$ and tangent to all surfaces displaced by a fixed distance from $\Gamma$. Then, the intrinsic divergence is equivalent to the standard divergence at the surface, $\nabla_\Gamma\cdot\boldsymbol{v} = \nabla\cdot\boldsymbol{v}$.
\end{principle}

Combining the two principles, higher order surface derivatives can be replaced with standard Cartesian derivatives in the embedding space \cite{marz2012calculus}. Thus, surface PDEs can be extended in the embedding space by replacing the surface derivatives with the corresponding standard ones using the two principles. Standard discretization techniques can be used for the solution of the PDE in the embedded space, while interpolation identifies the solution on the surface.

Standard finite differences \cite{ruuth2008simple} or finite differences based on radial basis functions (RBF-FD) \cite{petras2018rbf} have been used in the embedding space for the solution of PDEs on surfaces. However, instabilities have been observed in the application of these discretizations to stiff PDEs with large time step-size, which led to the introduction of different techniques of imposing the solution in the embedded space to be constant in the normal direction \cite{macdonald2009implicit,petras2019least,von2013embedded}. However, it is not clear how to choose the proper penalty term as described in \cite{von2013embedded} for stabilization of general stencils, and the least squares stabilization approach described in \cite{petras2019least} can be computationally expensive.

The idea of embedding the surface on a higher space dimension and using intrinsically a constant-along-normal extension in the embedding space was first presented in \cite{piret2012orthogonal} using global RBF interpolation. This is the so-called Radial Basis Functions Orthogonal Gradients method (RBF-OGr) and an improved sparse version based on RBF-FD was presented in \cite{piret2016fast}. Locally embedded stencils for solving  PDEs on folded surfaces were introduced in \cite{cheung2015localized}, which use intrinsically a constant-along-normal extension in the embedding space and the closest point principles to calculate the surface derivatives in a meshfree setting. Other meshfree works use embedding conditions \cite{chen2019meshless,cheung2018kernel} or an orthogonal projection operator \cite{chen2019meshless,chen2023kernel,chen2020extrinsic,tang2021localized} for kernel-based collocation methods for the solution of PDEs on surfaces. The main advantage of such methods is the independence of the grid in the embedded space, thus sparing the need of a closest point function, which can be quite challenging to compute. Gaussian RBFs were used therein to approximate the differential operators using RBF-FD.

In this work, we introduce local embedded stencils for the approximation of surface derivatives in combination with PHS+poly RBF-FD for the in-surface interpolation. Unlike \cite{cheung2015localized}, this work shows that it suffices to intrinsically extend along the normal direction only the reference node to overcome the rank deficiency of the polynomial basis, thus avoiding intersection of information along the normals when approximating high curvature areas. An extensive parameter analysis is presented regarding the size of the extension and the distance of the points placed in the embedded space along the normal direction. Our numerical results demonstrate the stability, flexibility and the ease in obtaining high order approximations, without the need for stabilization techniques. The method is applied on Turing pattern formation on rather general point clouds, as well as on advection PDEs. Extensions to PDEs on moving surfaces are also presented, with applications on strongly coupled Turing pattern systems.

The remainder of the paper is organized as follows. In Section \ref{sec:RBFFD}, we provide a short review on PHS+poly RBF-FD approximations. Section \ref{sec:method} provides an in-depth exposition of the proposed method, accompanied by a thorough numerical parameter exploration. In Section \ref{sec:numResults} we solve several surface PDE problems including forced diffusion, Turing patterns on either implicit surfaces or point clouds, and advection on the sphere and torus. This approach is extended to PDEs on moving surfaces in Section \ref{sec:movingSurfaces}. Section \ref{sec:conclusions} concludes the discussion with some final remarks.

\section{Review of PHS+poly generated RBF-FD approximations} \label{sec:RBFFD}

Consider a set of $N$ nodes $X=\{\boldsymbol{x}_i\}_{i=1}^N\subset\Omega\subseteq\mathbb{R}^d$ and a continuous target function $u:\Omega\rightarrow\mathbb{R}$ sampled at the nodes in $X$, i.e. $$u(\boldsymbol{x}_i)=u_i,\quad i=1,\dots,N.$$ 
Consider $\sigma^k\subset\{1,2,\dots,N\}$ as a subset comprising the indices of the $n$ closest points to a point $\boldsymbol{x}_k\in X$ (including itself), and let $\sigma^k_j$ represent the index of the $j$th closest point to $\boldsymbol{x}_k$ within this set. In finite difference (FD) approximations, the action of a linear differential operator $\mathcal{L}$ on the target function $u$ at $\boldsymbol{x}_k$ is approximated as a linear combination of the function values,
\begin{eqnarray}
\left. \mathcal{L}u \right|_{\boldsymbol{x}_k} \approx \sum_{j=1}^n w_j u(\boldsymbol{x}_{\sigma^k_j}),\label{RBFFD_approx}
\end{eqnarray}
where $w_j$ are the FD weights and $\boldsymbol{x}_{\sigma^k_j}$ is the $j$-th closest node to $\boldsymbol{x}_k$.

In radial basis functions-generated finite difference (RBF-FD) approximations, the weights $w_j$ are obtained 
enforcing \eqref{RBFFD_approx} to be exact for radial basis functions interpolants of the form 
\begin{equation}
s(\boldsymbol{x}) = \sum_{j=1}^n \lambda_j\,\phi(||\boldsymbol{x}-\boldsymbol{x}_{\sigma^k_j}||) + \sum_{i=1}^{L} \beta_i \, p_i(\boldsymbol{x}).
\label{RBF_int}
\end{equation}
where $\phi(r)$ is a RBF function centered at $\boldsymbol{x}_{\sigma^k_j}$, 
$\{p_i(\boldsymbol{x})\}_{i=1}^L$ is a monomial basis for the multivariate polynomial space $\Pi_l^d$ of total degree $l$ in $d$ dimensions, with $L=\binom{l+d}{l}<n$, and $\lambda_j$ and $\beta_i$ are scalar parameters. In this work, the RBF function $\phi$ is chosen to be the so-called polyharmonic spline (PHS) function of degree $m$, i.e. $\phi(r) = r^m$, where $m=2q+1$, $q\in\mathbb{N}$. 
Equation \eqref{RBF_int}, together with the orthogonality conditions
\begin{equation}
\sum_{j=1}^n\lambda_j p_i(\boldsymbol{x}_{\sigma^k_j})=0,\quad i=1,\ldots,L,
\label{RBF_orth}
\end{equation}
form the following system of equations determining the RBF-FD weights,
\begin{eqnarray}
\begin{bmatrix}
A & P \\[0.25ex]
P^T  & 0 \\[0.25ex]
\end{bmatrix}
\begin{bmatrix}
  \boldsymbol{w}\\
  \boldsymbol{\gamma}
\end{bmatrix} =
\begin{bmatrix}
  \left.\mathcal{L}\boldsymbol{\phi}\right|_{\boldsymbol{x}_k}\\
  \left.\mathcal{L}\boldsymbol{p}\right|_{\boldsymbol{x}_k}
\end{bmatrix}, 
\label{RBFFD_col}
\end{eqnarray}
where 
\begin{eqnarray*}
\begin{array}{lcl}
(A)_{ij}=\phi(||\boldsymbol{x}_{\sigma^k_i}-\boldsymbol{x}_{\sigma^k_j}||),&& i,j=1,\dots,n,\\[2ex]
(P)_{ij} = p_j(\boldsymbol{x}_{\sigma^k_i}), && i=1,\dots,n,\,\, j=1,\dots,L,\\[2ex]
(\left.\mathcal{L}\boldsymbol{\phi}\right|_{\boldsymbol{x}_k})_{i}=\left.\mathcal{L}\phi(||\boldsymbol{x}-\boldsymbol{x}_{\sigma^k_i}||)\right|_{\boldsymbol{x}_k}, && i=1,\dots,n,\\[2ex]
(\left.\mathcal{L}\boldsymbol{p}\right|_{\boldsymbol{x}_k})_{i}=\left.\mathcal{L}p_i(\boldsymbol{x})\right|_{\boldsymbol{x}_k},&& i=1,\dots,L,
\end{array}
\end{eqnarray*}
 and $\boldsymbol{\gamma}$ are the Lagrange multipliers of the underlying minimization problem \cite{B19}. It can be shown that $A$ is guaranteed to be positive definite on the subspace of vectors in $\mathbb{R}^n$ satisfying the orthogonality conditions in \eqref{RBF_orth} whenever $0 \leq q \leq l$ is chosen \cite{fasshauer2007meshfree}. Substitution of \eqref{RBFFD_col} in the approximation \eqref{RBFFD_approx} yields 
\begin{eqnarray}\label{eqn_RBFapprox}
\left. \mathcal{L}u \right|_{\boldsymbol{x}_k} \approx 
\begin{bmatrix}
  \left.\mathcal{L}\boldsymbol{\phi}\right|_{\boldsymbol{x}_k} \quad \left.\mathcal{L}\boldsymbol{p}\right|_{\boldsymbol{x}_k}
\end{bmatrix}^T 
\begin{bmatrix}
A & P \\[0.25ex]
P^T  & 0 \\[0.25ex]
\end{bmatrix}^{-1}
\begin{bmatrix}
  \boldsymbol{u}\\
  \boldsymbol{0}
\end{bmatrix}. 
\end{eqnarray}
The RBF-FD weights are assembled into a $N\times N$ sparse differentiation matrix with $n$ non-zero entries per row, which approximates the differential operator $\mathcal{L}$ on the point cloud $X$. Observe also that the calculation of RBF-FD weights is a pre-processing step that only needs to be performed once for any fixed $N$ and has a computational cost $O(n^3N)$. It also requires as pre-processing a nearest-neighbor search to find the stencils, which costs $O(N \log N )$ operations when using a $kd$-tree algorithm. For more detailed information on PHS+poly generated RBF-FD approximations, readers are referred to \cite{fornberg2015primer}.

\section{An RBF-FD constant-along-normal method}\label{sec:method}

\subsection{Description of the method}

The goal of this section is to describe the algorithm we propose to approximate a linear differential operator defined on a smooth closed surface $\Gamma$ of codimension $1$ embedded in $\mathbb{R}^{3}$ (which can be easily generalized to $\mathbb{R}^{d}$) with RBF-FD type approximations based on PHS, augmented with polynomials of degree $l$. 

\begin{figure}
\centering
\includegraphics[width=\textwidth]{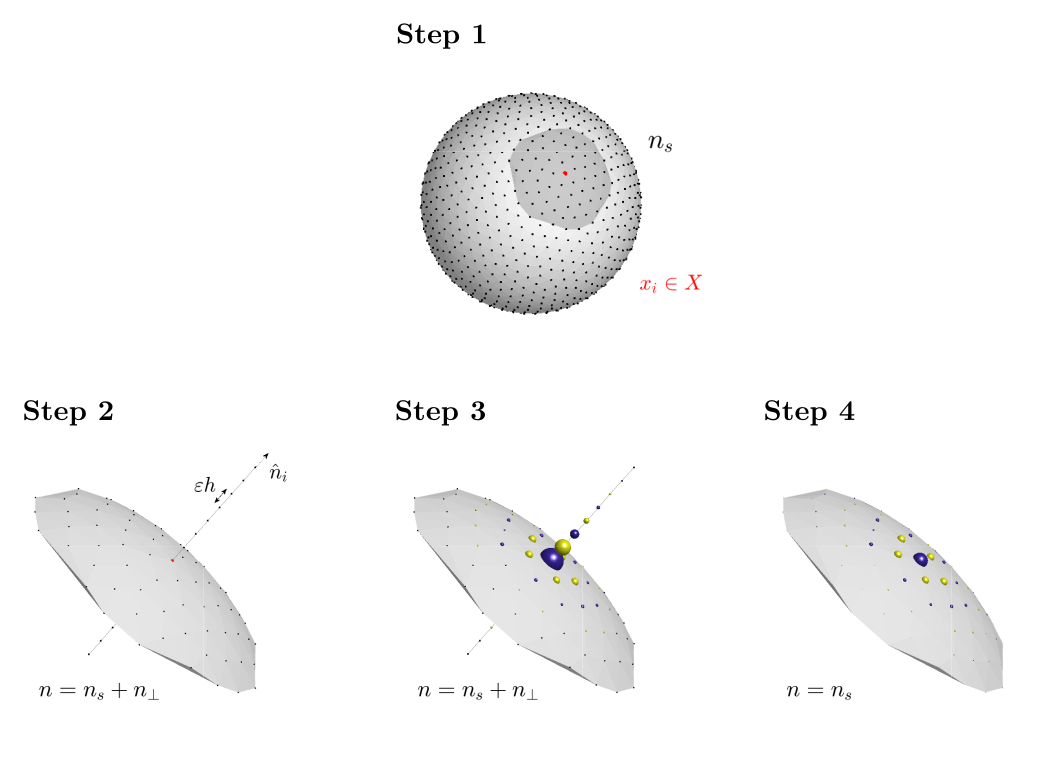}
\caption{A visualization of our method for calculating RBF-FD weights at a reference node $\boldsymbol{x}_i$ (in red). The introduced stencils use $n_s$ surface nodes (within the dark gray area) and $n_{\perp}$ out-of-surface nodes along the normal direction $\vec{n}_i$ to the reference node $\boldsymbol{x}_i$, with a distance of $\varepsilon h$ between each node. Note the exaggerated spacing between nodes along the normal direction in this figure (as $0<\varepsilon < 1$). The weights of the total $n$ nodes ($\boldsymbol{w}_{s}$ for the on-surface nodes and $\boldsymbol{w}_{\perp}$ for the out-of-surface nodes) in the stencil are shown as spheres, with blue indicating negative values and yellow indicating positive values. The radius of the spheres corresponds to the square root of the magnitude of the weights. Using the closest point principles, the final stencil ($\hat{\boldsymbol{w}}_s$) consists of only $n_s$ on-surface weights.}\label{fig:algorithm}
\end{figure}

For illustrative purposes, consider a linear differential operator defined on a smooth surface $\Gamma$, 
\begin{equation}\label{Lap}
\mathcal{L}_{\Gamma} u,
\end{equation}
where $u:\,\Gamma\to\mathbb{R}$ is a differentiable function. Let $\tilde{u}$ define a constant-along-normal extension of $u$ in the embedding space, such that $\tilde{u}|_\Gamma = u$. Using Principles \ref{p1} and \ref{p2}, a standard differential operator $\mathcal{L}$ applied on $\tilde{u}$ coincides with the action of the surface operator $\mathcal{L}_\Gamma$ on $u$ at the surface $\Gamma$, i.e.
\begin{equation}\label{embeddedLap}
\mathcal{L}_{\Gamma} u = \mathcal{L} \tilde{u}|_\Gamma.
\end{equation}
Let $X=\{\boldsymbol{x}_i\}_{i=1}^N$ denote the set of scattered node locations sampling the surface $\Gamma$, with unit normal vectors $\boldsymbol{\hat{n}}_i$ and average internodal distance $h$, and let $\left.u\right|_X$ denote the function values at $X$. For each ${\boldsymbol x}_i\in X$, which we refer to as \emph{reference node}, the proposed algorithm (visualized in Figure \ref{fig:algorithm}) follows the steps:

\vspace*{4ex}

\begin{enumerate}[leftmargin=2cm, labelindent=0cm, labelsep=5pt, font=\normalfont,itemsep=6ex]

\item[\bf Step 1.] Find the $n_s$ closest points 
		$$X_i=\{\boldsymbol{x}_j\}_{j\in\sigma_i}\subset X,\quad|\sigma_i|=n_s,$$
where $n_s=2\binom{l+d}{l}$ as suggested in \cite{BFFB17,shankar2020robust}, unless stated otherwise. Here $|\cdot|$ denotes the cardinality of the set.

\item[\bf Step 2.] Extend an equispaced node layout in the embedding space formed by $n_{\perp}$ (even) nodes in the normal direction of $\boldsymbol{x}_i$, 
\begin{eqnarray*}
X_i^{\perp}= \left\lbrace {\boldsymbol x}_i \pm j(\varepsilon h)\boldsymbol{\hat{n}}_i\right\rbrace_{j=1}^{n_{\perp}/2}, \qquad
\left\lbrace
\begin{array}{l}
n_{\perp}\geq l+1, \quad l \mbox{ is odd},\\
n_{\perp}> l+1, \quad l \mbox{ is even},
\end{array}
\right.
\end{eqnarray*}
where $h$ is the average internodal distance between the nodes in $X$ and $0<\varepsilon<1$ is a small scaling parameter. Recall that $l$ is the degree of the augmented polynomial term in the RBF interpolant \eqref{RBF_int}. 

\item[\bf Step 3.] Use $X_i\cup X_i^{\perp}$ as stencil to compute the RBF-FD weights as in \eqref{eqn_RBFapprox}, such that
			\begin{eqnarray}
				\left.\mathcal{L}\tilde{u}\right|_{{\boldsymbol x}_i} \approx
				\begin{bmatrix}				
				\boldsymbol{w}_s^T & \boldsymbol{w}_{\perp}^T				
				\end{bmatrix}
				\begin{bmatrix}
				\left.\tilde{u}\right|_{X_i} \\ 
				\left.\tilde{u}\right|_{X_i^{\perp}}
				\end{bmatrix},\label{step3}
			\end{eqnarray}
where $\boldsymbol{w}_s$ denotes the weights at $X_i$ and $\boldsymbol{w}_{\perp}$ at $X_i^{\perp}$. Observe that $n=|X_i\cup X_i^{\perp}|=n_s+n_{\perp}$.

\item[\bf Step 4.] By construction $\tilde{u}|_{X_i} = u|_{X_i}$. Since $\tilde{u}$ is constant-along-normal extension of $u$ in the embedding space, it follows that $\tilde{u}|_{X_i^{\perp}} = u|_{\boldsymbol{x}_i}$. Therefore, the RBF-FD approximation in equation \eqref{step3} can be rewritten as 
		\begin{eqnarray*}
		\left.\mathcal{L}\tilde{u}\right|_{{\boldsymbol x}_i} 
		&\approx&
		\hat{\boldsymbol{w}}_s^T \left.u\right|_{X_i} 
		\end{eqnarray*}
		where
\begin{eqnarray*}
(\hat{\boldsymbol{w}}_s)_{_k}=
\left\lbrace
\begin{array}{ll}
(\boldsymbol{w}_s)_{_1} + \sum_{j=1}^{n_{\perp}} (\boldsymbol{w}_{\perp})_{_j}, &\quad \mbox{if } k = 1,\\[4ex]
(\boldsymbol{w}_s)_{_k}, & \quad \mbox{if }k=2,\dots,n_s,
\end{array}
\right.\\
\end{eqnarray*}
and the stencil size reduces to $n=n_s$.
\end{enumerate}

The local weights $\hat{\boldsymbol{w}}_s$ for $\boldsymbol{x}_i$ are assembled into the $i$-th row of the sparse differentiation matrix, which contains only $n_s$ nonzero entries. The steps are repeated for all nodes on the surface.

Observe that this approach only requires nodes at scattered node locations on the surface and its corresponding normals. The constant-along-normal extension reduces the stencil size from $n_s+n_{\perp}$ collocation nodes to $n_s$ differentiation weights in step (4), leading to $n_s$ non-zero elements per row in the differentiation matrix. It does not require the calculation of the closest points, and can be used to approximate any linear differential operator in any space dimension over smooth surfaces defined either implicitly/parametrically or by point clouds. The nodes along the normal extension are placed using an average internodal distance $h$ on $X$, scaled by $0<\varepsilon<1$. The number $n_{\perp}$ depends on the polynomial degree $l$. At least, it is required that $n_{\perp}\geq l+1$. For a given polynomial degree $l$, the convergence of the approach when approximating a linear differential operator $\mathcal{L}$ of order $k$ is $\mathcal{O}(h^{l+1-k})$, provided that the interpolation approximation is $ \mathcal{O}(h^{l+1})$ \cite{B19,BFFB17,FFBB16}.

Note that surfaces with significant curvature changes over small patches can cause off-surface points to inadvertently penetrate the surface. While techniques to mitigate this issue may involve varying the scaling parameter \(\varepsilon\) for the normal direction extension or selecting stencil nodes with consistent normals, as discussed in \cite{petras2022meshfree}, such methods are beyond the scope of this work. Nonetheless, the next section demonstrates that reducing \(\varepsilon\) has minimal impact on the accuracy of the proposed method. And not less important, if the surface is discretized sufficiently fine to capture its features, there is no need to adjust parameters, as the scaling of out-of-surface nodes relative to the internodal distance on the surface will be inherently accurate.

\subsection{Parameter exploration}

The proposed method introduces two new parameters that influence the constant-along-normal extension in the embedding space: the number of out-of-surface nodes, $n_{\perp}$, and the distance between them, $\varepsilon h$, where $h$ is the average internodal distance. To assess their impact and clarify their effect on the performance of our method, we designed the following test cases: the approximation of the Laplace-Beltrami operator on the unit circle in $\mathbb{R}^2$, the approximation of the surface Poisson's equation on the "tooth" model \cite{lehto2017radial} and the unit sphere in $\mathbb{R}^3$, and the approximation of the Laplace-Beltrami operator on high-curvature surfaces. 


\subsubsection{A 2-D test case to illustrate the effect on the contourlines}

\begin{figure}
\centering

\includegraphics[width=\textwidth,trim={10ex 0ex 10ex 0ex},clip]{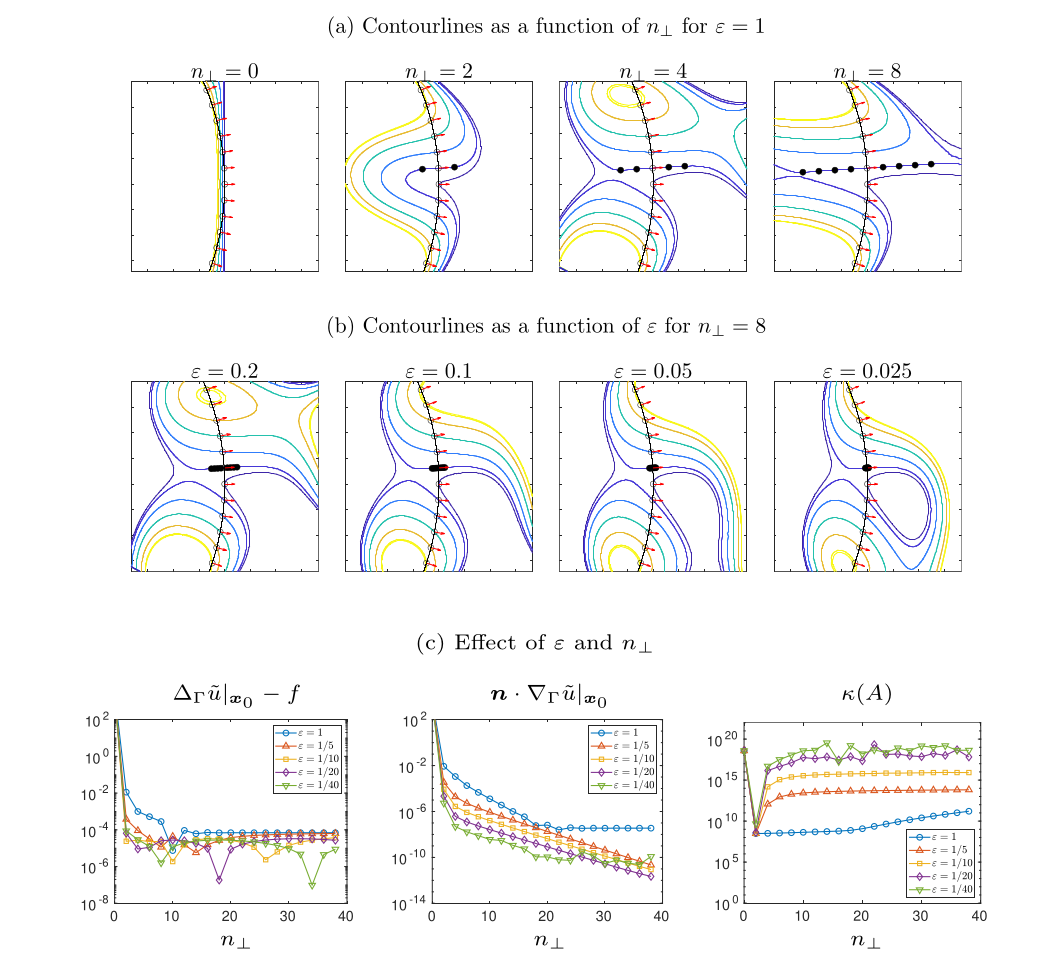}

\caption{Effect of the constant extension along the normal when approximating the surface Laplacian of \(u(x,y) = \sin(\pi x)\cos(\pi y)\) on the unit circle \(\Gamma\) at \((0.9980, 0.0628)\) with \(n_s = 20\) using PHS \(r^7\) with polynomials of degree 2: (a) Contour lines vs. \(n_{\perp}\) for \(\varepsilon=1\); (b) Contour lines vs. \(\varepsilon\) for \(n_{\perp}=8\); (c) Effect of varying \(\varepsilon\) and \(n_{\perp}\) on the residual of the surface Laplacian, projection of surface gradient on normal direction and condition number of the collocation system.}\label{fig:param2D}
\end{figure}

To illustrate the algorithm, we have approximated the surface Laplacian of the test function $u(x,y)=\sin(\pi x)\cos(\pi y)$ 
at a point $\boldsymbol{x}_0$ on the unit circle $\Gamma$.  Figure \ref{fig:param2D} shows the isocontours for varying parameters $n_{\perp}$ and $\varepsilon$, using $n_s=20$ and PHS $r^7$ augmented with polynomials of degree $l=2$. Observe in panel (a) how the contourlines become parallel to the normal at the stencil center as $n_{\perp}$ increases for a fixed $\varepsilon=1$; or, in panel (b), how the contourlines close around the stencil center for a fixed $n_{\perp}=8$ and $\varepsilon=0.2,\,0.1,\,0.05$ and 0.025. 

Panel (c) illustrates three key aspects: the residual of the approximation $\Delta_{\Gamma}\tilde{u}|_{\boldsymbol{x}0}-f$ (first row), the projection of the surface gradient onto the normal direction $\boldsymbol{\hat{n}}\cdot\nabla_{\Gamma} \tilde{u}|_{\boldsymbol{x}0}$ (second row), and the condition number of the RBF-FD collocation matrix $\kappa(A)$ (third row). It is important to note that the method becomes ill-conditioned when $n_{\perp}=0$ (i.e., no equispaced node layout along the normal direction), due to the linear dependence of the polynomial terms ${1,x^2,y^2}$, which results in $x^2+y^2=1$. This issue is avoided when $n_{\perp}\geq l+1$. Observe also that by either increasing $n_{\perp}$ or decreasing $\varepsilon$, the term $\boldsymbol{\hat{n}}\cdot\nabla_{\Gamma} \tilde{u}|_{\boldsymbol{x}0}$ decreases. This is expected because $\boldsymbol{\hat{n}}$ and $\nabla_{\Gamma} \tilde{u}|_{\boldsymbol{x}0}$ become orthogonal in this limit, as shown in panels (a) and (b). However, while the condition number increases, it does not appear to significantly impact the residual of the approximation, provided that $n_{\perp}\geq l+1$.


\subsubsection{Test case in 3-D}\label{sec_param3Dcase}

\begin{figure}
\centering
\hspace*{-4ex} \includegraphics[scale=0.4,trim={4ex 2ex 2ex 2ex},clip,rotate=90]{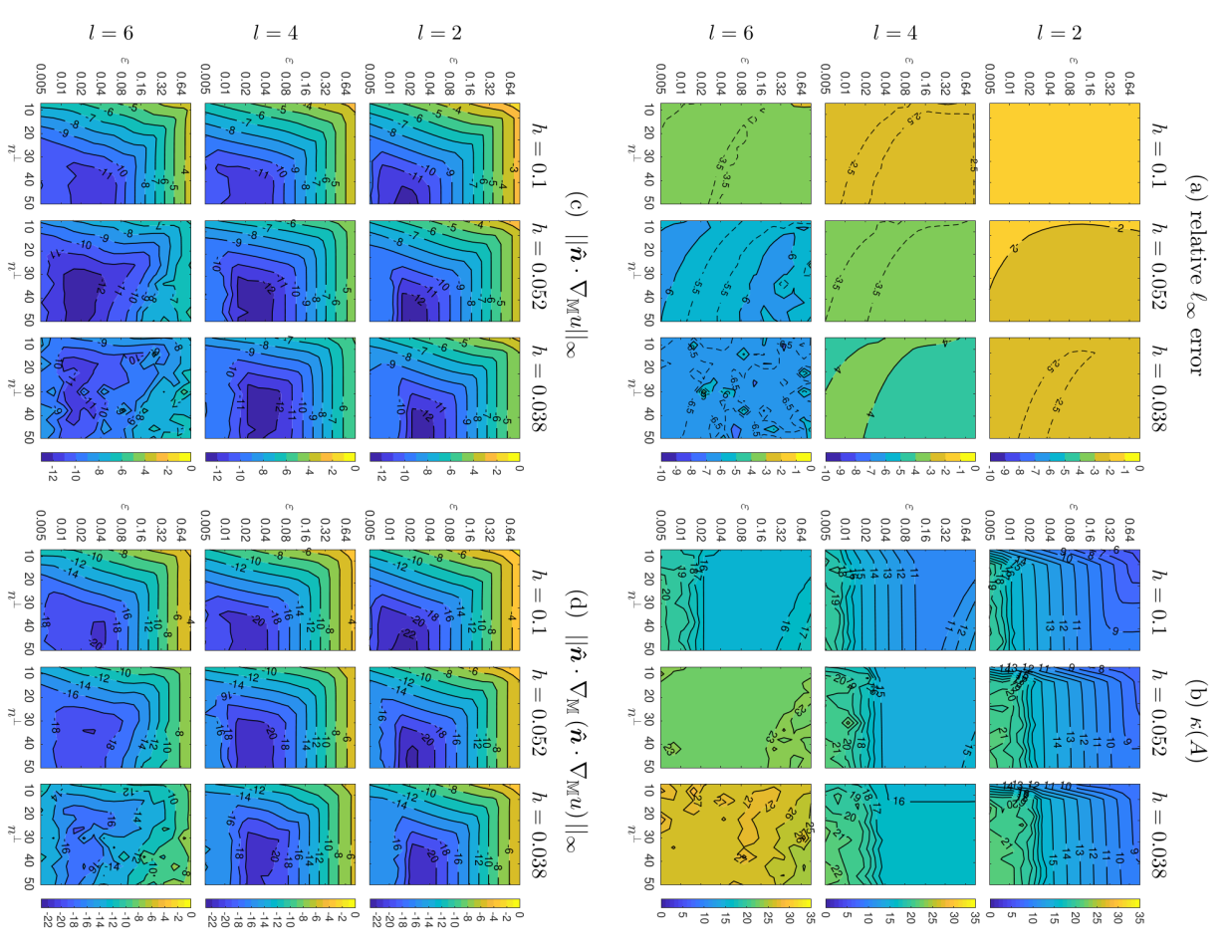}
\caption{Effect of varying $\varepsilon$ against $n_{\perp}$ when solving \eqref{BVPPDE} using the constant-along-normal extension with PHS $r^5$ augmented with polynomials of degree $l$ on the (a) accuracy, (b) condition number, (c) projection of the surface gradient on the normal direction and (d) $||\boldsymbol{\hat{n}}\cdot\nabla_{\Gamma} \left( \boldsymbol{\hat{n}}\cdot\nabla_{\Gamma} u \right)||_{\infty}$ for different resolutions (rows) and augmented polynomial degrees $l$ (columns). The numbers on the colorbars represent powers of 10.}
\label{fig:param3D}
\end{figure}

Consider the PDE problem
\begin{equation}
-\Delta_{\Gamma} u = f\label{steadyPDE}
\end{equation}
defined on a closed surface of codimension $1$ embedded in $\mathbb{R}^{3}$, and $u:\,\Gamma\to\mathbb{R}$ is a differentiable function. 
To explore the parameter dependence of our approach, we have considered a synthetic BVP test where the boundary effect (stencil deformation near boundary) is prevented. This is given by 
\begin{eqnarray}
\begin{array}{rcll}
-\Delta_{\Gamma} u &=& f, \quad &z\geq0 \\
u &=& g, \quad &z<0
\end{array}\label{BVPPDE}
\end{eqnarray}
where the smooth surface is the so-known tooth model \cite{lehto2017radial},
\begin{eqnarray}
\Gamma = \{(x,y,z)\in\mathbb{R}^{3} \,\,|\,\, x^8 + y^8 + z^8 - (x^2 + y^2 + z^2) = 0 \},\label{tooth}
\end{eqnarray}
and $f$ and $g$ are computed exactly from the test function
\begin{eqnarray}
u_1(x,y,z) = -\cos\left(\frac{3}{4}\pi x\right)\cos(\pi y)\sin\left(\frac{3}{2}\pi z\right).\label{testfun_u1}
\end{eqnarray}
In this problem, we have used PHS $r^5$ augmented with polynomials of degree $l$ and varied the different parameters involved in the method. For instance, Figure \ref{fig:param3D} shows the effect of varying $\varepsilon$ against $n_{\perp}$ on (a) the relative $\ell_{\infty}$ error, (b) the condition number $\kappa(A)$, (c) the projection $||\boldsymbol{\hat{n}}\cdot\nabla_{\Gamma} u ||_{\infty}$ and (d) $||\boldsymbol{\hat{n}}\cdot\nabla_{\Gamma} \left( \boldsymbol{\hat{n}}\cdot\nabla_{\Gamma} u \right)||_{\infty}$. Each of these quantities has been computed with different augmented polynomial degrees $l=2,4,6$ (rows) and resolutions $h=0.1,0.076,0.052,0.036$ (columns) on the domain $(\varepsilon,n_\perp) \in [0.005,0.9] \times [9,50]$. Observe the following:

\begin{itemize}

\item Figure \ref{fig:param3D}-(a): the effect of $\varepsilon$ and $n_{\perp}$ in the accuracy is almost negligible. The decrease of $h$ or the increase of $l$ has a stronger impact.

\item Figure \ref{fig:param3D}-(b): the condition number of the collocation matrix increases as $\varepsilon$ decreases, while $n_{\perp}$ does not seem to have any effect as long as $n_{\perp}\geq l+1$. Apparently, the high condition number computed only affects $\boldsymbol{\hat{n}}\cdot\nabla_{\Gamma} \left( \boldsymbol{\hat{n}}\cdot\nabla_{\Gamma} u\right)$ for $\varepsilon<0.02$ and any $l$. The effect on the accuracy seems to be negligible except for $l=6$. This might reflect a limitation in the definition of the condition number, being sensitive to scaling issues as noted in \cite{FFBB16}.


\item Figure \ref{fig:param3D}-(c)\&(d): the dependence of $\boldsymbol{\hat{n}}\cdot\nabla_{\Gamma} u$ and $\boldsymbol{\hat{n}}\cdot\nabla_{\Gamma} \left( \boldsymbol{\hat{n}}\cdot\nabla_{\Gamma} u\right)$ on $\varepsilon$ and $n_{\perp}$ is very similar. Both decrease as $\varepsilon$ decreases until the ill-conditioned region for $\varepsilon<0.02$ is reached, and display a saturation region for $\varepsilon$ low and $n_{\perp}$ large.

\end{itemize}

\begin{figure}
\centering
\includegraphics[scale=0.33,trim={0ex 2ex 0ex 2ex},clip]{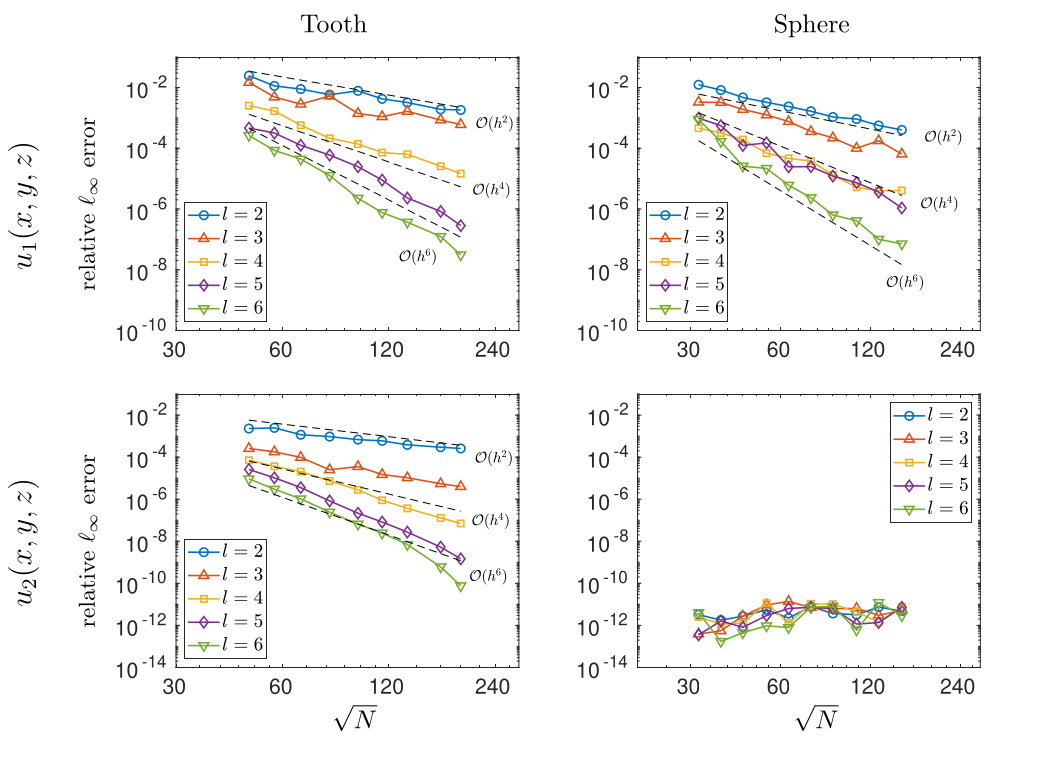} 
\caption{Convergence as a function of $\sqrt{N}$ when using the constant-along-normal extension to solve the BVP \eqref{BVPPDE} on the tooth \eqref{tooth} and unit sphere for the test functions $u_1$ \eqref{testfun_u1} and $u_2(x,y,z)=-xy$ with $n_{\perp}=10$ and $\varepsilon=0.05$.}
\label{Fig:conv_ToothSphere}
\end{figure}

\vspace*{2ex}

Additionally, we computed the convergence rates when solving problem \eqref{BVPPDE} as a function of $\sqrt{N}$ for $r^5$ augmented with different polynomial degrees $l$, fixing $n_{\perp}=10$ and $\varepsilon=0.05$. Figure \ref{Fig:conv_ToothSphere} shows the corresponding results on the tooth model \eqref{tooth} (first column) and on the unit sphere (second column) for the same test function $u_1$ (first row). Observe that in both surfaces the convergence rate is determined by the polynomial degree $l$ and improves as $l$ increases (as expected for PHS+poly approximations \cite{B19,BFFB17,FFBB16}). In particular, $l=2,3$ seems to achieve $\mathcal{O}(h^2)$; $l=4,5$ yields $\mathcal{O}(h^4)$; and $l=6$ gets $\mathcal{O}(h^6)$. 

In the second row, we have computed the problem for the polynomial test function
\begin{eqnarray*}
u_2(x,y,z)=-xy.
\end{eqnarray*}
On the tooth model, the convergence rates keep the same as for $u_1$, but the accuracy improves (the curves shift down). On the sphere, our approach achieves machine precision. This is related with the fact that we are using $r^5$ augmented with polynomials of degree $l\geq2$ to approximate a second degree polynomial on the implicit surface $x^2+y^2=1$. Further results on the spectra of the Laplace--Beltrami operator $\Delta_{\Gamma}$ as the polynomial degree increases are considered in Appendix \ref{App:Leig}.

\newpage


\subsubsection{Approximation on High-Curvature Surfaces}

In this section, we analyze the performance of our approach on high-curvature surfaces, where a surface feature may be very close to another part of the surface along the normal direction. In particular, we aim to understand the trade-off between the out-of-surface nodes along the normal direction and the proximity of different parts of the surface in this direction.

\begin{figure}
\centering

\includegraphics[scale=0.33,trim={4ex 2ex 4ex 2ex},clip]{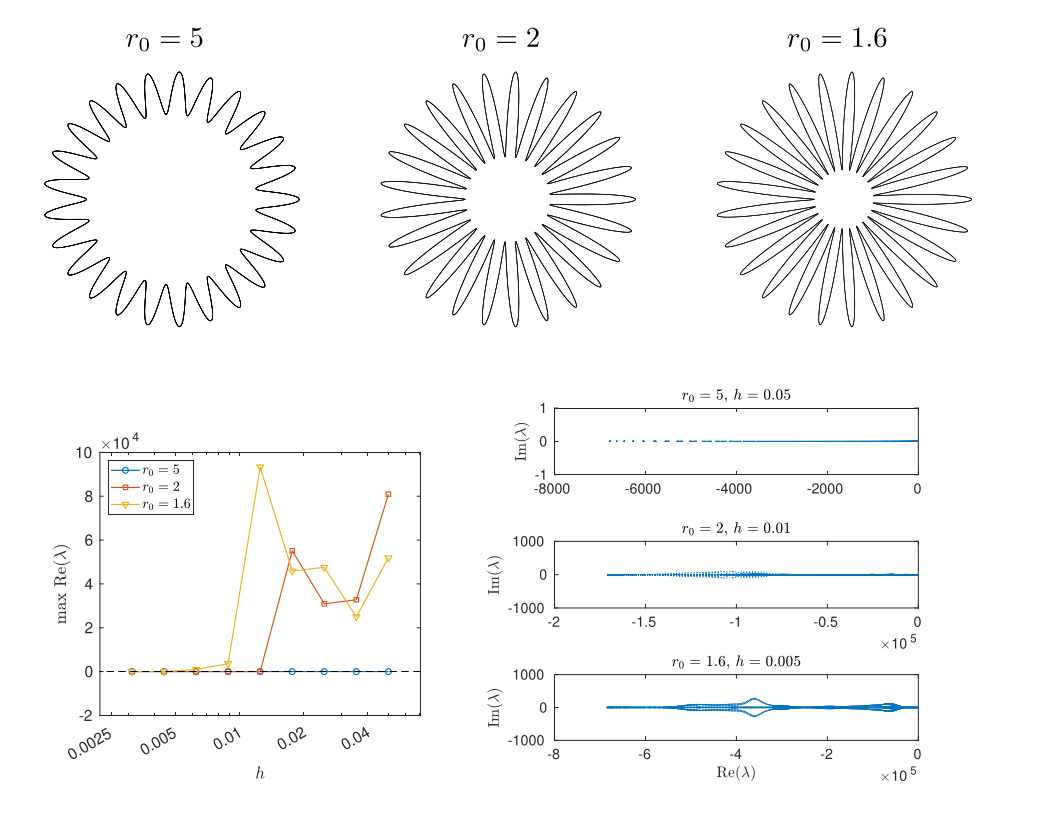}

\caption{Illustration of the rose curve for different values of the radial offset \( r_0 \) (top row) and the corresponding eigenvalue spectra of the numerical approximation of the Laplace-Beltrami operator for different internodal distances \( h \) (bottom row). The curves are defined by \( k = 25 \) and \( r_0 = 5, 2, \), and \( 1.6 \).}
\label{Fig:2D_hCurv}
\end{figure}

For our first numerical test, we approximate the Laplace-Beltrami operator on the so-called \textit{rose curve}, expressed in polar coordinates as
\[
(x(\theta), y(\theta)) = \left( r_0 + \cos(k \theta) \right) (\cos(\theta), \sin(\theta)),
\]
where \( r_0 \) and \( k \) are parameters that control the pattern of the rose curve. Specifically, \( r_0 \) is the radial offset that shifts the rose curve outward from the origin, while \( k \) determines the number of petals. In this test, we consider \( k = 25 \) for \( r_0 = 5, 2, \) and \( 1.6 \), as depicted in the first row of Figure \ref{Fig:2D_hCurv}. Note that as \( r_0 \) approaches 1, regions of higher curvature appear with features very close to each other along the normal direction. We have discretized these curves by computing the arc length of the curve and setting nodes to a fixed internodal distance \( h \). For each discretization, we have approximated the Laplace-Beltrami operator using PHS of degree 3 augmented with polynomials of degree 2, with a stencil size of \( n_s = 9 \) and fixed parameters \( \varepsilon = 0.1 \) and \( n_{\perp} = 4 \). We analyze the eigenvalue spectra. The left plot in the second row of Figure \ref{Fig:2D_hCurv} shows the $\max \mbox{Re}(\lambda)$ against the internodal distances used \( h \). Observe that as \( h \) decreases and the curves are discretized finely enough, the differentiation matrix provides a good approximation of the surface Laplacian, with none of the eigenvalues crossing into the right half-plane (i.e. $\max \mbox{Re}(\lambda) < 0$). We have depicted the corresponding eigenvalue spectra in these cases for \( r_0 = 5, h = 0.05 \); \( r_0 = 2, h = 0.01 \); and \( r_0 = 1.6, h = 0.005 \) in the second plot of the figure. Note that the imaginary part of all eigenvalues approaches zero.

This suggests that as long as the curve is finely discretized with respect to its features (high-curvature regions), the parameters \((\varepsilon, n_{\perp})\) do not need to be tuned. The key seems to be properly scaling the out-of-surface nodes along the normal direction relative to the internodal distance on the surface, i.e.,
\[X_i^{\perp} = \left\lbrace {\boldsymbol x}_i \pm j(\varepsilon h)\boldsymbol{\hat{n}}_i \right\rbrace_{j=1}^{n_{\perp}/2}.\]
In this case, \( X_i^{\perp} \) is simply formed by a layer of 4 nodes at distances \( 0.1h \) and \( 0.2h \) from the reference point. Thus, once the regions of the surface with high curvature are properly resolved, \( X_i^{\perp} \) avoids having off-surface points inadvertently penetrating the surface.

\begin{figure}
\centering
\includegraphics[scale=0.33,trim={10ex 0ex 10ex 0ex},clip]{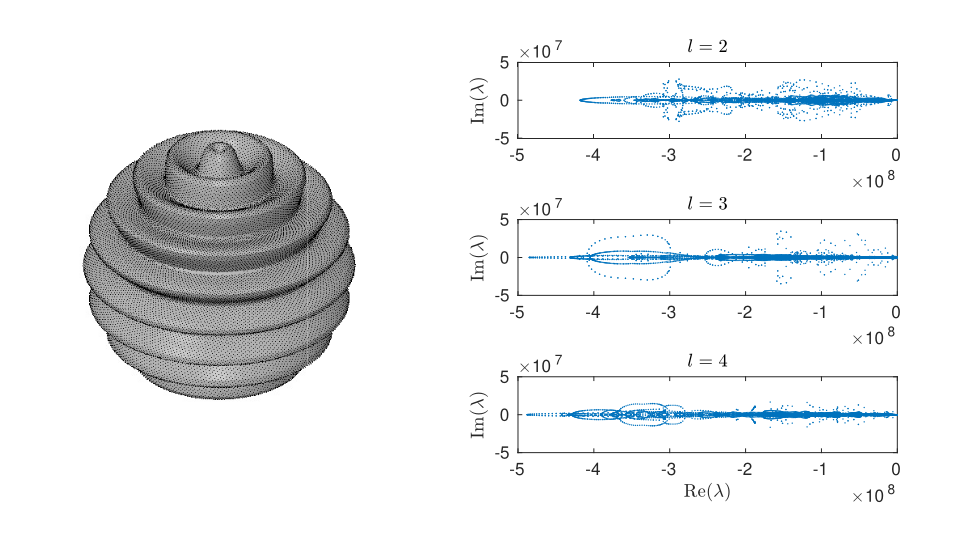}
\caption{Left: Visualization of the bumpy sphere. Right: Eigenvalue spectra of the numerical approximation of the Laplace-Beltrami operator for the bumpy sphere, using PHS of degree 3 augmented with polynomials of degree \( l = 2, 3, \) and 4, a stencil size \( n_s = \left\lfloor{\frac{3}{2}\binom{l+3}{3}}\right\rfloor \), and fixed parameters \(\varepsilon = 0.1\) and \( n_{\perp} = 11 \).}
\label{Fig:3D_hCurv}
\end{figure}

Following the same idea, we consider the bumpy sphere shown in Figure \ref{Fig:3D_hCurv}, which is a sphere whose radius is parametrized by \( r(\phi) = 1 + \gamma \sin(k \phi) \), where \(\gamma = 0.1\) and \(k = 21\). On a discretization formed by \(N = 50{,}000\) nodes, we approximate the Laplace-Beltrami operator using PHS of degree 3 augmented with polynomials of degree \( l = 2, 3, \), and 4, with a stencil size \( n_s = \left\lfloor{\frac{3}{2}\binom{l+3}{3}}\right\rfloor \) and fixed parameters \(\varepsilon = 0.1\) and \(n_{\perp} = 11\). The eigenvalue spectra are displayed on the right side of Figure \ref{Fig:3D_hCurv}. As in the previous 2D case, the eigenvalue spectra demonstrate that with a fine enough discretization of the surface, a good approximation of the Laplace-Beltrami operator can be achieved without the need for tuning the parameters \((\varepsilon, n_{\perp})\).

\section{Numerical results}\label{sec:numResults}
\subsection{Surface diffusion}

In this section we analyze the convergence of our numerical method for approximating the heat equation with forcing 
\begin{equation}\label{forcedHeatPDE}
u_t = \Delta_\Gamma u + f(t,u),
\end{equation}
on two surfaces: the unit sphere and the torus. In all the tests, we choose $n_\perp=14$ and $\varepsilon=0.2$, and the stencil size $n_s=\left\lfloor{\frac{3}{2}\binom{l+3}{3}}\right\rfloor$. We use PHS $r^5$ augmented with polynomials of degree $l$.

\begin{figure}
\centering
\includegraphics[trim={0ex 0ex 5ex 3ex},clip,scale=0.33]{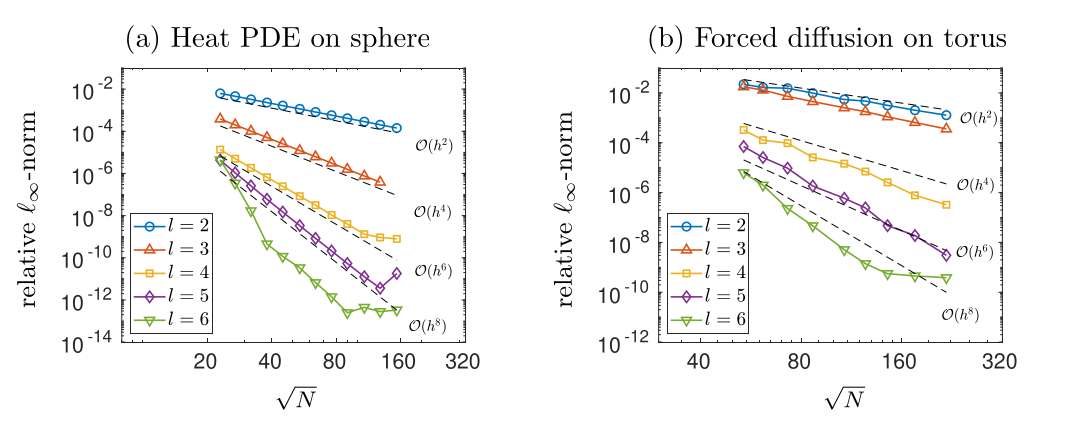}
\caption{The convergence of the method when solving \eqref{forcedHeatPDE} on two different surfaces in $\mathbb{R}^3$ with $n_\perp=14$ and $\varepsilon=0.2$.}\label{Fig:surfDiff}
\end{figure}

\subsubsection{Diffusion on the Sphere}

To begin, the equation \eqref{forcedHeatPDE} is defined on the surface of the unit sphere in $\mathbb{R}^3$. The exact solution is given as a series of spherical harmonics $$u(t,\theta,\phi) = \frac{20}{3\pi} \sum_{l=1}^{\infty} e^{-l^2/9} e^{-t\,l(l+1)}Y_{ll}(\theta,\phi),$$ where $\theta$ and $\phi$ are longitude and latitude, respectively, and $Y_{lm}$ is the degree $l$ order $m$ spherical harmonic. Here, the forcing term is $f(t,u)=0$. Since the coefficients decay rapidly, the series is truncated after 30 terms as in \cite{macdonald2009implicit}. We have evolved the numerical solution in time until $t=0.5$ using RK4 with a time-step $\Delta t=0.5/N$. 

Figure \ref{Fig:surfDiff} (left) shows the convergence results when combining $r^5$ with polynomials of degree $l=2,3,\dots,6$. Observe that convergece improves as the polynomial degree increases. Actually, the error decreases until the machine precision is reached and starts diverging as $\mathcal{O}(h^2)$.

\subsubsection{Forced Diffusion on a Torus}

This test is similar to the one presented in \cite{fuselier2013high}, where the forcing function $f(t,u)$ in \eqref{forcedHeatPDE} is determined analytically from the exact solution
$$u(t,\boldsymbol{x}) = \frac{1}{8} e^{-5t} x \left(x^4 - 10x^2y^2 + 5y^4\right)\left(x^2 + y^2 -60z^2\right),$$
where $\boldsymbol{x}=(x,y,z)$ is a point on the Torus. As in the previous case, we have evolved the numerical solution in time until $t=0.2$ using RK4 with a time-step $\Delta t=0.5/N$. The numerical error is analyzed as a function of the spatial resolution for different augmented polynomial degrees. Figure \ref{Fig:surfDiff} (right) shows the convergence results for different augmented polynomial degrees.

\subsection{Reaction-diffusion systems: Turing patterns}\label{sec_TuringStatic}

\begin{figure}[h!]
\centering

\includegraphics[trim={4ex 5ex 4ex 5ex},clip,scale=0.38,rotate=90]{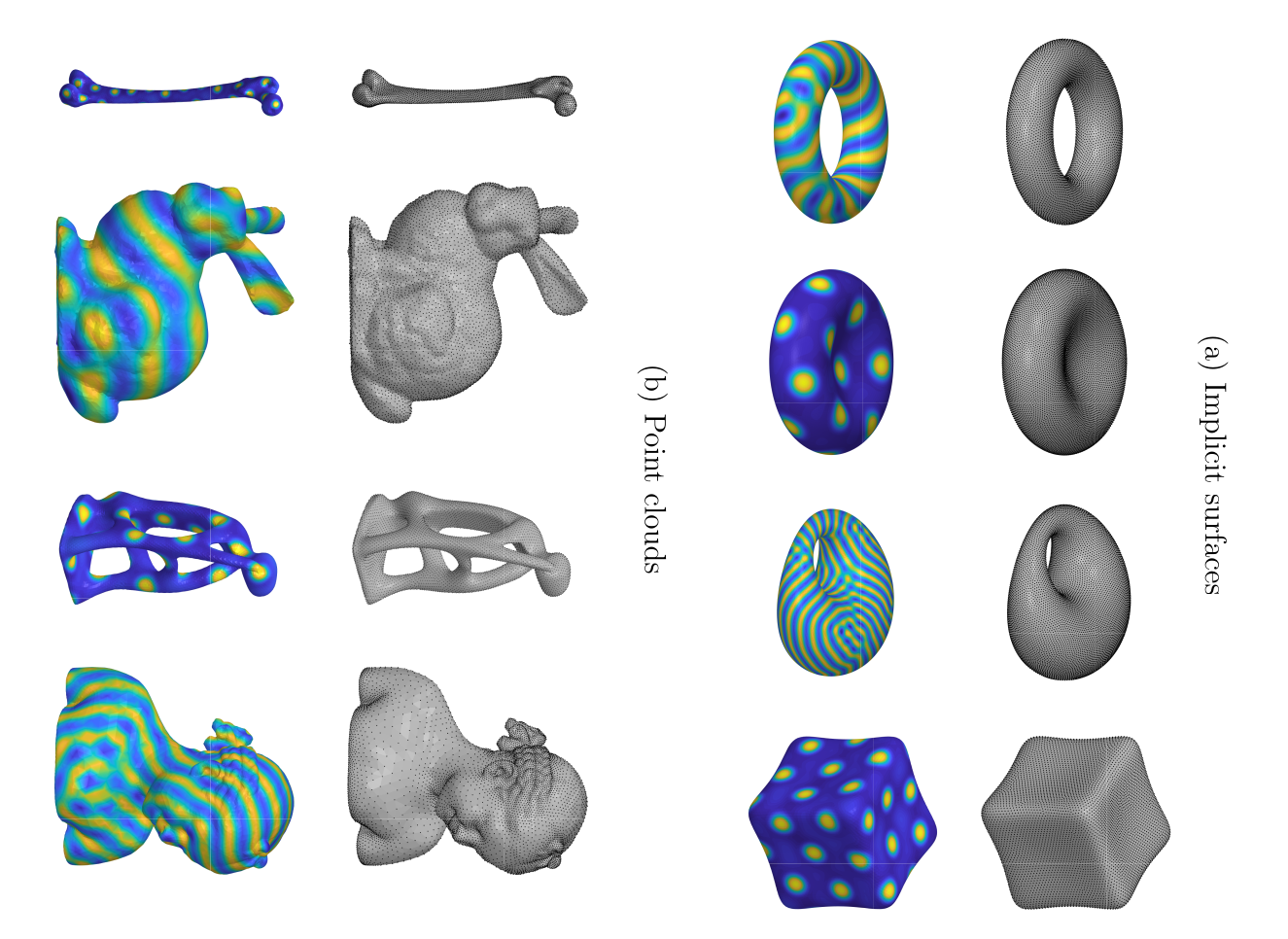}

\caption{Steady Turing spot and stripe patterns resulting from solving equations \eqref{Turing1} and \eqref{Turing2} on: (a) Implicit surfaces (torus, red blood cell, Duphin's cyclide and tooth models with $9{,}458$, $13{,}778$, $11{,}884$ and $9{,}632$ points); (b) Point clouds (Femur, Stanford Bunny, Woodthinker and Bimba models with $3{,}954$, $9{,}948$, $8{,}260$ and $12{,}391$ points).}
\label{Fig:Turing}
\end{figure}

We now present applications of the proposed method to solving reaction-diffusion equations on more general surfaces. Following the numerical tests in \cite{fuselier2013high,lehto2017radial,shankar2015radial}, we solve the two-species Turing system
\begin{eqnarray}
\frac{\partial u}{\partial t} &=& \delta_u\Delta_{\Gamma}u + \alpha u (1-\tau_1 v^2) + v(1-\tau_2u),\label{Turing1}\\
\frac{\partial v}{\partial t} &=& \delta_v\Delta_{\Gamma}v + \beta v (1+\frac{\alpha\tau_1}{\beta}uv) + u(\gamma+\tau_2v),\label{Turing2}
\end{eqnarray}
on surfaces defined implicitly by algebraic expressions, such as the torus, red blood cell model, Dupin cyclide and tooth model; and on more general surfaces without an implicit/parametric representation, such as the Femur, Stanford Bunny, Woodthinker and Bimba models obtained using JIGSAW \footnote{JIGSAW is a Delaunay-based unstructured mesh generator for two- and three-dimensional geometries. The models considered in these work are freely available at \texttt{https://github.com/dengwirda/jigsaw-models}.} \cite{engwirda2016off}. The results obtained using the parameters from Table~\ref{tbl_params} are shown in Figure~\ref{Fig:Turing}. Light and dark shading indicates high and low concentrations, respectively.

In these examples, the Laplace--Beltrami operator is approximated using PHS $r^5$ augmented with polynomials of degree 6, and parameters $n_{\perp}=10$ and $\varepsilon=0.1$. For the implicit surfaces, the normals $\boldsymbol{\hat{n}}$ are computed exactly. For the point clouds, we approximate them using \texttt{pcnormals} from the Computer Vision Toolbox in MATLAB. To advance in time, we use one step of the first order Semi--implicit Backward Difference Formula (SBDF), followed by SBDF2 with a time-step $\Delta t=0.02$. These are custom functions implemented in MATLAB. At each time step, the resulting implicit system is solved using MATLAB's built-in BiCGSTAB solver with a tolerance $10^{-11}$ and incomplete LU with 0 level of fill in as preconditioner. The corresponding eigenvalue distribution of the Laplace--Beltrami operator is attached in the Appendix \ref{App:Teig}. Observe that our method provides stable approximations on both implicit surfaces and point clouds, with all eigenvalues having a negative real part for all surfaces considered.

\begin{table}
\centering
\caption{The parameters of the Turing system used to obtain spots and stripes in our results with $\delta_u=0.516\delta_v$}\label{tbl_params}
\begin{tabular}{cccccccc}
\hline
\textbf{Pattern} & $\delta_v$ & $\alpha$ & $\beta$ & $\gamma$ & $\tau_1$ & $\tau_2$ & \textbf{Final time}\\
\hline
Spots & 4.5e-3 & 0.899 & -0.91 & -0.899 & 0.02 & 0.2 & 600\\
Stripes & 2.1e-3 & 0.899 & -0.91 & -0.899 & 3.5 & 0 & 6000\\
\hline
\end{tabular}
\end{table}

\subsection{Surface advection} 

In this section, we investigate the use of our RBF-FD method on the surface advection equation 
\begin{equation}
\begin{array}{l}\displaystyle
\frac{\partial q(\boldsymbol{x},t)}{\partial t} + \boldsymbol{v}\cdot\nabla_{\Gamma} \, q(\boldsymbol{x},t)=0
\end{array}\label{surfAdv}
\end{equation}
which represents the transport of a tracer $q(\boldsymbol{x},t)$ with initial condition $q(\boldsymbol{x},0)=q_0(\boldsymbol{x})$ over a surface $\Gamma$ by a given velocity field $\boldsymbol{v}$. We implement two examples from the bibliography: the solid-body rotation over the sphere \cite{gunderman2020transport} and the transport over the torus in a time-independent, spatially varying flow field \cite{shankar2020robust}.

In all cases, the velocity field should drive the solution back to the initial condition after a given time period. Thus, the exact solution at the final time is known, allowing us to compute errors. We analyze the convergence rates 
on quasi-uniform node sets of increasing sizes using PHS $r^3$ augmented with polynomials of degree $l=3,4,5$ and $6$. We use maximal determinant (MD) nodes on the sphere, as reported in \cite{gunderman2020transport}, and staggered nodes on the torus as in \cite{shankar2020robust}. 

\begin{figure}
\centering

\includegraphics[scale=0.35,trim={4ex 0ex 4ex 0ex},clip]{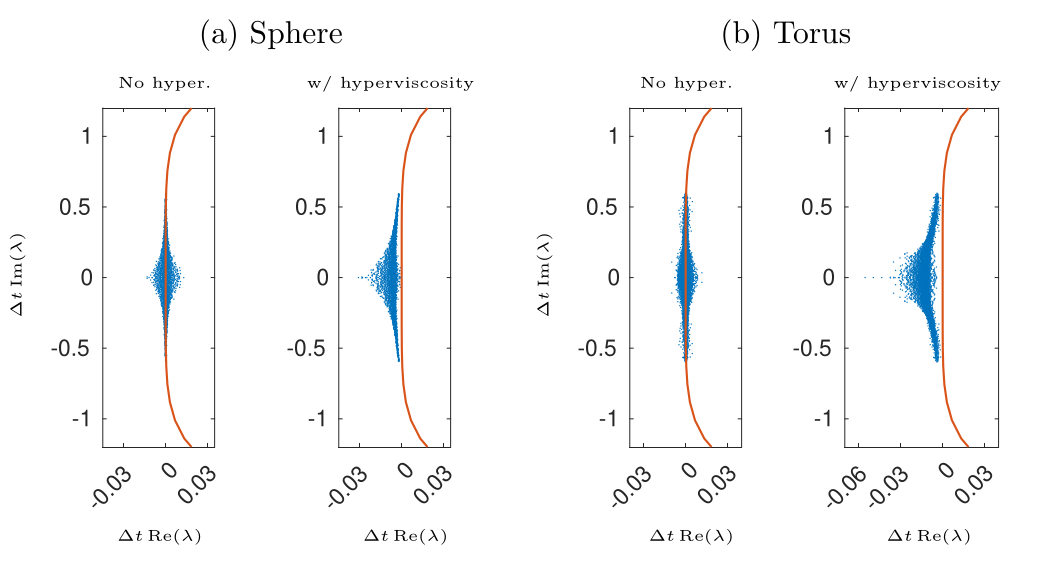} 

\caption{The effect of hyperviscosity on the eigenvalues of the spatial differentiation matrix approximating $\boldsymbol{v}\cdot\nabla_{\Gamma}$ with $\boldsymbol{v}=\boldsymbol{1}$, on (a) the Sphere with $N = 2{,}916$ MD nodes, and (b) the Torus with $N = 5{,}312$ staggered nodes. The line in red represents the RK4 stability region.}\label{Fig:eigenAdv}
\end{figure}

It is well-known that RBF-FD differentiation matrices, which approximate advective terms, usually contain spurious eigenvalues with small positive real parts that can cause instabilities when time-stepping, especially when an explicit time discretization is used. To overcome this issue, we follow the same approach as in \\cite{fornberg2011stabilization,flyer2012guide,shankar2020robust} and add an artificial hyperviscosity term to the right-hand side of equation \eqref{surfAdv}, leading to
\begin{equation}
\begin{array}{l}\displaystyle
\frac{\partial q(\boldsymbol{x},t)}{\partial t} + \boldsymbol{v}\cdot\nabla_{\Gamma} \, q(\boldsymbol{x},t)=\gamma_k \Delta_{\Gamma}^k q(\boldsymbol{x},t),
\end{array}\label{surfAdv2}
\end{equation}
where $k= \lfloor\log(n_s)\rfloor$, $\Delta_{\Gamma}$ is the surface Laplacian and $\gamma_k = \epsilon\,(-1)^{k+1}h^{2k+1}$. Here, $h\sim1/\sqrt{N}$ is the average internodal distance of the node distribution and $\epsilon$ is a scaling factor which requires to be found by trial-and-error. In the present test cases, we found that $\epsilon=0.01$ for the torus and $\epsilon=0.001$ for the sphere works well. As a result of this little artificial hyperviscosity term, all spurious eigenvalues of the original differentiation matrix move into the negative half eigenplane and we can fit all them within the RK4 time stability region as shown in Figure \ref{Fig:eigenAdv}. The time step $\Delta t$ is chosen as stated below.

\subsubsection{Advection on the sphere}

As a first test case, we analyze the solid body rotation \cite{gunderman2020transport}, which involves the advection of an initial condition in a steady velocity field. The components of the steady velocity field for this test case in spherical coordinates are given by
\begin{eqnarray}
u(\lambda, \theta) = \sin(\theta) \cos(\lambda)\sin(\alpha ) - \cos(\theta) \cos(\alpha ),\quad v(\lambda, \theta) = \cos(\lambda) \sin(\alpha),\label{sph_field}
\end{eqnarray}
where $-\pi\leq\lambda\leq\pi$, $-\pi/2\leq  \theta\leq\pi/2$ and $\alpha$ is the angle of rotation with respect to the equator. We choose $\alpha=\pi/2$, which corresponds to advecting the initial condition over the poles. The flow field returns the solution to its initial position after a time $T = 2\pi$.

We test two different initial conditions. The first one, a compactly supported cosine bell centered at $\boldsymbol{p}=(1,0,0)$,
\begin{eqnarray}
q(\boldsymbol{x},0) =
\left\lbrace
\begin{array}{ll}
\frac{1}{2}\left(1+\cos\left(\frac{\pi r}{R_b}\right)\right) & r<R_b\\
0 & r\geq R_b 
\end{array}\right.\nonumber
\end{eqnarray}
where $r=\arccos(x)$ and $R_b=1/3$. This initial condition is $\mathcal{C}^1(\mathbb{S}^2)$ as it has a jump in the second derivative. The second initial condition is a Gaussian bell centered at the same point $\boldsymbol{p}$ and given by
\begin{equation}
q_2(\boldsymbol{x},0) = e ^{-6||\boldsymbol{x}-\boldsymbol{p}||^2}.
\end{equation}
This initial condition is $\mathcal{C}^{\infty}(\mathbb{S}^2)$, which makes it ideal to test high-order convergence. In both cases, we run one full revolution over the sphere, which corresponds to $T = 2\pi$,
using a time step $\Delta t = \frac{T}{10\sqrt{N}}$.

\begin{figure}
\centering
\includegraphics[trim={0ex 4ex 0ex 0ex},clip,scale=0.33]{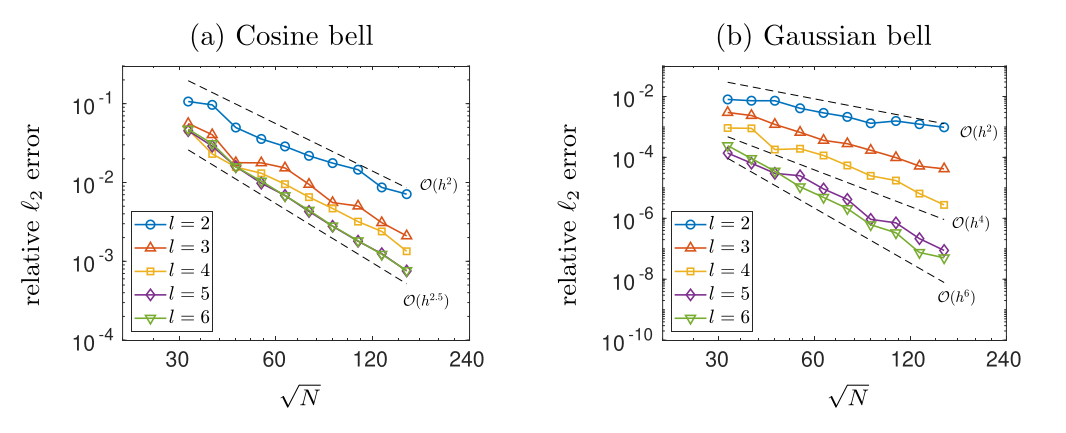}
\caption{Convergence on the sphere for the surface advection equation of a cosine bell (left) and a Gaussian bell (right) in a steady flow. The figure displays the relative $\ell_2$ error as a
function of $\sqrt{N}$ for polynomial degrees $l=2,3,4,5$ and $6$. The dashed lines indicate different convergence rates.
}\label{Fig:advSphere}
\end{figure}

Figure \ref{Fig:advSphere} shows the convergence on the sphere as a function of $\sqrt{N}$ for the surface advection equation of a cosine bell (left) and a Gaussian bell (right) in the steady flow \eqref{sph_field}. Observe that for the $\mathcal{C}^1(\mathbb{S}^2)$ cosine bells, increasing the polynomial degree $l$ does not increase convergence rates but does improve the accuracy, as has been observed among others in \cite{gunderman2020transport}. This is because it has a jump in the second derivative. For the $\mathcal{C}^{\infty}(\mathbb{S}^2)$ Gaussian bells, it is shown that increasing the polynomial degree $l$ results in a higher convergence rate, which is in closer agreement with the predicted values. In all cases, the stencil size chosen is $n_s=\left\lfloor{\frac{3}{2}\binom{l+3}{3}}\right\rfloor$, the number of points along the normal is $n_\perp=14$ and the scaling parameter $\varepsilon=0.2$. Observe that our method leads to competitive results if compared with those reported in \cite{gunderman2020transport} for $l=5$ and $n_s=37$.

\subsubsection{Advection on the torus}

In this section, we consider the experiment reported in \cite{shankar2020robust}, where a field $q(\boldsymbol{x},t)$ is advected over a torus of inner radius $1/3$ and outer radius 1 by the time-independent velocity field $\boldsymbol{v}=(u,v,w)$,
\begin{eqnarray}
\begin{array}{lcl}
u &=& -\frac{2}{3}\sin(2(\phi-\theta))\cos(3\phi) - 3\left[1 + \frac{1}{3}\cos(2(\phi-\theta))\right]\sin(3\phi),\\[3ex]
v &=& -\frac{2}{3}\sin(2(\phi-\theta))\sin(3\phi) + 3\left[1 + \frac{1}{3}\cos(2(\phi-\theta))\right]\cos(3\phi),\\[3ex]
w &=& -\frac{2}{3}\cos(2(\phi-\theta)),\\
\end{array}\label{torus_field}
\end{eqnarray}
expressed with respect to the Cartesian basis. This velocity field will advect a particle on the surface of the torus around a $(3,2)$ torus knot, returning to its initial position after a time $T=2\pi$. 

As in the previous experiment, we consider here two different initial conditions. The first corresponds to a pair of compactly supported cosine bells centered at $\boldsymbol{p}_1 = (1+1/3, 0, 0)$ and $\boldsymbol{p}_2=-\boldsymbol{p}_1$, given by
$$q(\boldsymbol{x},0) = 0.1 + 0.9(q_1 + q_2 ),$$
where $q_1$ and $q_2$ are defined as
\begin{eqnarray}
q_{1,2}=\left\lbrace
\begin{array}{ll}
\frac{1}{2}\,(1 + \cos(2\pi r_{1,2} (\boldsymbol{x}))) & \mbox{if} \quad r_{1,2}(\boldsymbol{x}) < 0.5,\\[2ex]
0 & \mbox{if} \quad r_{1,2}(\boldsymbol{x}) \geq 0.5,
\end{array}
\right.
\end{eqnarray}
and $r_{1,2}(\boldsymbol{x}) = \left\| \boldsymbol{x}-\boldsymbol{p}_{1,2} \right\|$. This solution has only one continuous derivative.

The second initial condition consists of a pair of Gaussian bells centered at the same points $\boldsymbol{p}_{1,2}$ and given by the initial condition  
\[q(x,y,z,0) = e^{-a(x+(1+1/3)^2) +y^2 -1.5az^2} + e^{-a(x-(1+1/3))^2+y^2-1.5az^2},\]
where $a=20$. This is an infinitely smooth solution which allows us to test high-order convergence as a function of the augmented polynomial degree $l$. We set the time step to $\Delta t = \frac{T}{10\,v_{max} \sqrt{N}}$. Here, $v_{max}=4.1\approx\max_{\theta,\phi}\|\boldsymbol{v}\|_{\infty}$ is roughly the maximum magnitude of the velocity field over the torus. In all cases, we choose $n_s=\left\lfloor2.2\binom{l+3}{3}\right\rfloor$, $n_\perp=14$ and $\varepsilon=0.5$.
We select this stencil size because of some numerical instabilities at high resolutions for high polynomial degrees $l=5, 6$. Increasing it to this value effectively mitigated these instabilities.

\begin{figure}
\centering
\includegraphics[trim={0ex 4ex 0ex 0ex},clip,scale=0.33]{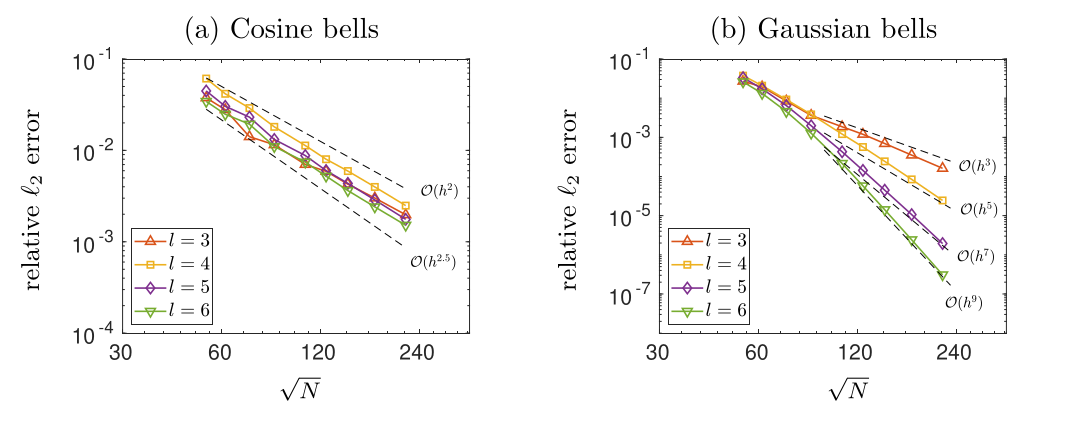} 
\caption{Convergence on the torus for the surface advection equation of two cosine bells (left) and two Gaussian bells (right) in a steady flow. The figure displays the relative $\ell_2$ error as a function of $\sqrt{N}$ for polynomial degrees $l=3,4,5$ and $6$. The dashed lines indicate different convergence rates.}
\label{Fig:advTorus}
\end{figure}

Figure \ref{Fig:advTorus} shows the convergence on the torus as a function of $\sqrt{N}$ for the surface advection equation of two cosine bells (left) and two Gaussian bells (right) in the steady flow \eqref{torus_field}. As in the sphere, increasing the polynomial degree $l$ does not increase convergence rates for the cosine bells, but does improve the accuracy. However, it controls the convergence rate for the infinitely smooth Gaussian bells. Again, our numerical approach leads to very competitive results if compare with the results reported on \cite{shankar2020robust}.

\section{Extension to PDEs on moving surfaces}\label{sec:movingSurfaces}

In this section, we provide an algorithm that allows the extension of our method for the solution of PDEs on moving surfaces. Our numerical results include a convergence study of the method on diffusion conservation laws and Turing patterns on deforming surfaces.

\subsection{Algorithm description}

For the numerical solution of PDEs on moving surfaces, we perform a simple coupling with a particle movement algorithm. Higher order coupling methods require further analysis, which goes beyond the scope of this paper. 
The proposed algorithm uses the following steps:\\[2ex]


\begin{enumerate}[leftmargin=2cm, labelindent=0cm, labelsep=5pt, font=\normalfont,itemsep=2ex]

\item[\bf Step 1.] \textbf{Solution of PDE.} Solve the PDE using the method provided in Section~\ref{sec:method}.

\item[\bf Step 2.] \textbf{Particle motion.} Move the particles along the normal direction according to a motion law. The solution of the PDE is carried along with the particles, as constant-along-normal. 

\item[\bf Step 3.] \textbf{Resampling of particles.} Typical particle movement algorithms provide a resampling of the particles with addition/removal of points.

\item[\bf Step 4.] \textbf{In-surface interpolation.} Use any surface interpolation method to obtain the solution on the resampled surface.\\[2ex]

\end{enumerate}

The algorithmic Steps 2 and 3 are carried out using any particle tracking algorithm, such as \cite{chen2020kernel,leung2009grid,petras2016pdes,suchde2019fully}. In this paper, we consider the grid based particle method for the surface movement \cite{leung2009grid}.

The algorithmic Step 4 uses an in-surface interpolation technique for the identification of the solution on the updated surface points. Any meshfree interpolant could be considered, including \cite{cheung2015localized,shankar2018rbf,Suchde2019}. In our work, it is natural to consider the interpolant that is derived from Section~\ref{sec:method} by replacing the differential operator with the identity operator.

\subsection{Numerical results}
In our numerical results on evolving surfaces, unless stated otherwise, the same parameters as in the static surface examples in Section~\ref{sec:numResults} are considered. Since the point clouds on moving surfaces may change at every time step with an arbitrary fill distance, the points on the stencil are chosen with a minimum distance of $\Delta x/2$ among each other, similar to a previous work \cite{petras2022meshfree}.

\subsubsection{Diffusion conservation law on an expanding sphere}
To test the convergence of the method, we use a benchmark example as presented in \cite{suchde2019fully}. Consider the conservation law of a scalar quantity $u$ with a diffusive flux on a surface $\Gamma$ given by
\begin{equation}\label{eqn_consLaw}
\frac{Du}{Dt} + u\nabla_\Gamma\cdot\boldsymbol{v} - \Delta_\Gamma u = f,
\end{equation}
where $D/Dt$ is the total derivative, $\boldsymbol{v}$ is the surface velocity and $f$ is a forcing term. Let $\Gamma$ be the unit sphere, and $\boldsymbol{v}$ a velocity of the form 
$$\boldsymbol{v} = 0.5\boldsymbol{\hat{n}},$$
where $\boldsymbol{\hat{n}}$ is the outward unit normal vector. This motion leads to the expansion of the sphere, with a linear increase of its radius along time, i.e. $r(t) = 1+0.5t$. Using the function
$$u(t,x,y,z) = e^{-6t}xy,$$
as the desired solution of the conservation law \eqref{eqn_consLaw}, the forcing term is identified as
$$f(t,x,y,z) = \left(-6+\frac{2}{r(t)}+\frac{6}{r^2(t)}\right)u(t,x,y,z).$$
For the given constant-along-normal velocity, equation \eqref{eqn_consLaw} can also be written as
\begin{equation}\label{eqn_consLaw2}
\frac{\partial u}{\partial t} + 0.5\frac{\partial u}{\partial \boldsymbol{\hat{n}}} + 0.5\kappa u - \Delta_\Gamma u = f,
\end{equation}
where $\partial/\partial\boldsymbol{\hat{n}}$ is the normal derivative and $\kappa$ is the mean curvature. 

Let $X^0$ be the initial pointcloud of the unit sphere of size $N_0$, and $U^0$ the vector of the initial solution with $U^0_i=[u(0,x_i,y_i,z_i)]$, $\boldsymbol{x}_i = (x_i,y_i,z_i)\in X^0$, $i=1,...,N_0$. We numerically solve equation \eqref{eqn_consLaw2} using a first order implicit-explicit scheme \cite{ascher1995implicit}
$$(I-\Delta t \Delta_h)U^{n+1} = U^n + \Delta t(0.5K^n\circ U^n + F^n),$$
where $I$ is the identity matrix, $\Delta_h$ is the discretized Laplace-Beltrami operator as described in Section \ref{sec:method}, $U^{n+1}$ and $U^n$ is the vector of the numerical approximation of the solution at times $t_{n+1}$ and $t_n$ respectively, $K^n$ is the mean curvature at time $t_n$, $\circ$ is the elementwise product of two vectors and $F^n_i = [f(t_n,x_i,y_i,z_i)]$, $i=1,\dots,N_n$, with $N_n$ the size of the pointcloud at time $t_n$. Note that the method is Section \ref{sec:method} provides a constant-along-normal extension by construction, and additionally through Step 2 of our proposed algorithm the solution is carried along the normal direction in the particle motion. Thus, by construction $\partial u/\partial \boldsymbol{\hat{n}}=0$.

To show the spatial convergence of our method, we consider an evolution up to $t=0.5$ using a small time step-size of $\Delta t = 0.4\Delta x^2$, where $\Delta x$ is the grid size in the grid based particle method. Table~\ref{tbl:expandingSphere} shows the $\ell_\infty$-norm relative error of the numerical solution against the exact solution $u$ at the final time. 

\begin{table}
    \centering
    \caption{The convergence of the method for the diffusion conservation law on the expanding sphere using the $\ell_\infty$-norm relative error. Also, the initial $N_0$ and final (at $t=0.5$) $N_{0.5}$ number of points on the surface.}
    \begin{tabular}{lcclc}
        \hline
        $\Delta x$ & $N_0$ & $N_{0.5}$ & $\ell_\infty$-error & e.o.c. \\ \hline
         0.4 & 200 & 272 & 0.34 & - \\
         0.2 & 688 & 1088 & 0.086 & 1.99 \\
         0.1 & 2816 & 4400 & 0.022 & 1.96\\
         0.05 & 10976 & 17312 & 0.0057 & 1.97\\
        \hline
    \end{tabular}
    \label{tbl:expandingSphere}
\end{table}

\subsubsection{Turing patterns on evolving surfaces}
In our next example, we consider a surface cross-diffusion reaction-diffusion system as in \cite{petras2019least} of the form
\begin{eqnarray}\label{eqn_Turingmoving}
\frac{Du}{Dt} & = & \Delta_\Gamma u + d_w\Delta_\Gamma w - u\nabla_\Gamma\cdot\boldsymbol{v} + f_1(u,w), \nonumber\\[0.1cm]
& &\\[-0.5cm]
\frac{Dw}{Dt} & = & d_u\Delta_\Gamma u + c\Delta_\Gamma w - w\nabla_\Gamma\cdot\boldsymbol{v} + f_2(u,w), \nonumber
\end{eqnarray}
where $c$, $d_w$ and $d_u$ are positive constants and $f_1$ and $f_2$ are the mixing functions. Using a solution $u$ and mean curvature $\kappa$ dependent velocity along the normal direction of the form
$$\boldsymbol{v} = (-0.01\kappa+0.4u)\boldsymbol{\hat{n}}=: V(\kappa,u)\boldsymbol{\hat{n}},$$
and the analysis in \cite{Dziuk2013}, the system of equations \eqref{eqn_Turingmoving} can be rewritten as
\begin{eqnarray}\label{eqn_Turingmoving2}
\frac{\partial u}{\partial t} & = & \Delta_\Gamma u + d_w\Delta_\Gamma w - V(\kappa,u)u\kappa - V(\kappa,u)\frac{\partial u}{\partial \boldsymbol{\hat{n}}} + f_1(u,w), \nonumber\\[0.1cm]
& &\\[-0.5cm]
\frac{\partial w}{\partial t} & = & d_u\Delta_\Gamma u + c\Delta_\Gamma w - V(\kappa,u)w\kappa - V(\kappa,u)\frac{\partial w}{\partial \boldsymbol{\hat{n}}} + f_2(u,w). \nonumber
\end{eqnarray}
Similar to the previous example, we solve numerically the system \eqref{eqn_Turingmoving2} using a first order implicit-explicit scheme, written as a block system
\begin{eqnarray*}
\left(\begin{array}{cc}
    I -\Delta t\Delta_h & -d_w\Delta t\Delta_h \\
    -d_u\Delta t\Delta_h & I - c\Delta t\Delta_h
\end{array}\right)\left(\begin{array}{c}
     U^{n+1}  \\
     W^{n+1} 
\end{array}\right) = \left(\begin{array}{c}
      U^n + \Delta t(-V(K^n,U^n)\circ U^n\circ K^n + f_1(U^n,W^n)) \\
      W^n + \Delta t(-V(K^n,U^n)\circ W^n\circ K^n + f_2(U^n,W^n))
\end{array}\right),
\end{eqnarray*}
where the quantities are consistently defined as in the previous example. 

\begin{figure}
\centering
\includegraphics[width=\textwidth]{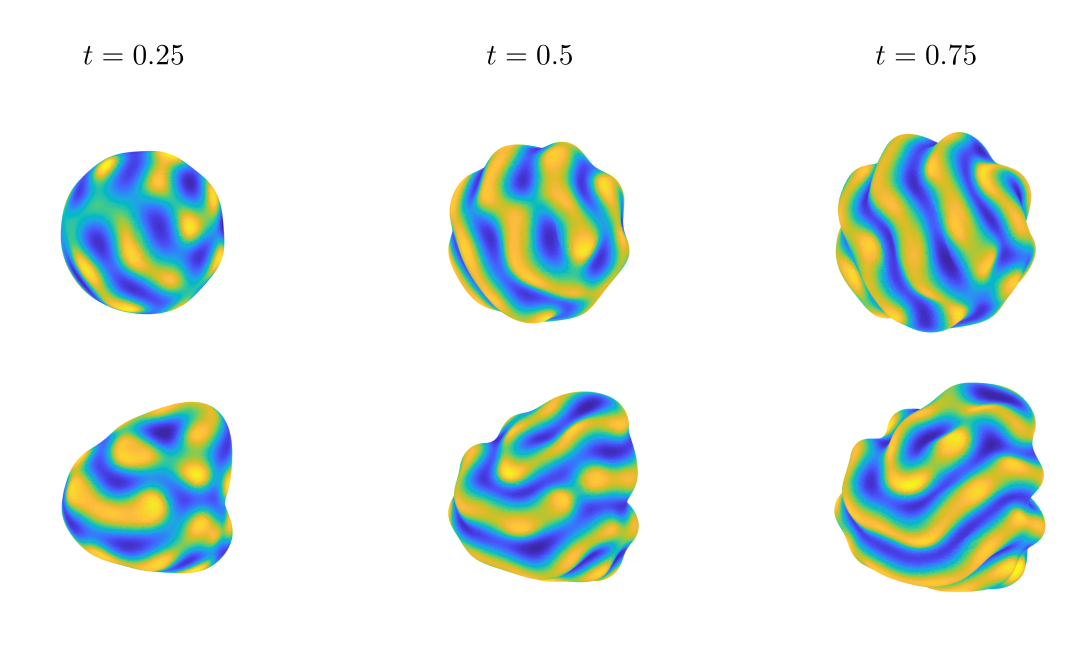}
\caption{Snapshots of the solution of the cross-diffusion reaction-diffusion system on the unit sphere (top) and the Dziuk's surface (bottom) at times $t=0.25$, $0.5$ and $0.75$.}
\label{Fig:movingSurfaces}
\end{figure}

In this example, we use the parameters $d_u=d_w=1$, $c=10$ and the mixing functions
$$f_1(u,w) = 200(0.1-u+u^2w)\quad\text{and}\quad f_2(u,w) = 200(0.9-u^2w).$$
Using as initial condition random perturbations around 1 and 0.9 for $u$ and $w$ respectively, we solve the system with $\Delta x=0.07$ and $\Delta t=0.002$ up to a final time of $t=0.75$. The chosen initial surfaces are the unit sphere and the Dziuk surface, defined as an implicit surface
$$(x-z^2)^2+y^2+z^2=1.$$
Snapshots of the solutions at different times appear in Fig.~\ref{Fig:movingSurfaces}.

\section{Conclusions} \label{sec:conclusions}
This paper presents a stable meshfree PHS+poly RBF-FD constant-along-normal method for solving PDEs on surfaces of codimension 1 embedded in $\mathbb{R}^3$ (which can be easily generalized for surfaces in $\mathbb{R}^{d}$). The method is built upon the principles of the closest point method, without the use of a grid or a closest point mapping. The key result of this approach is that  stable spatial discretizations for hyperbolic or stiff differential operators (such as the surface Laplacian) can be explicitly computed using stencils locally embedded in $\mathbb{R}^d$ at each point of the surface. Specifically, we demonstrate that it is enough to consider a constant extension along the normal direction only at the reference node to overcome the rank deficiency of the polynomial basis. 
The robustness of this approach is tested for a number of examples, including the advection and diffusion equations, as well as Turing pattern formations on surfaces defined implicitly or by point clouds. A simple coupling approach with a particle tracking method shows the potential of the proposed method for the solution of PDEs on evolving surfaces in the normal direction. The method appears robust in the different fill distances from the evolving point clouds.

\section*{Acknowledgements}
The work of VB has received support from the Madrid Government (Comunidad de Madrid-Spain) under the Multiannual Agreement with UC3M (H2SAFE-CM-UC3M). AP has been partially supported by the State of Upper Austria. SJR has been partially supported by the Natural Sciences and Engineering Research Council of Canada (RGPIN 2022-03302).  The authors thank Grady Wright for sharing with us the node distributions for the RBC, Dupin cyclide and "tooth" models used in the paper, and Varun Shankar for his help with the advection on the torus test case.

\newpage
\appendix

\section{Eigenvalue spectra for the Laplace--Beltrami operator}

\subsection{Parameter exploration: tooth and unit sphere} \label{App:Leig}

\begin{figure}
\centering
\includegraphics[scale=0.35,trim={0ex 0ex 0 0},clip]{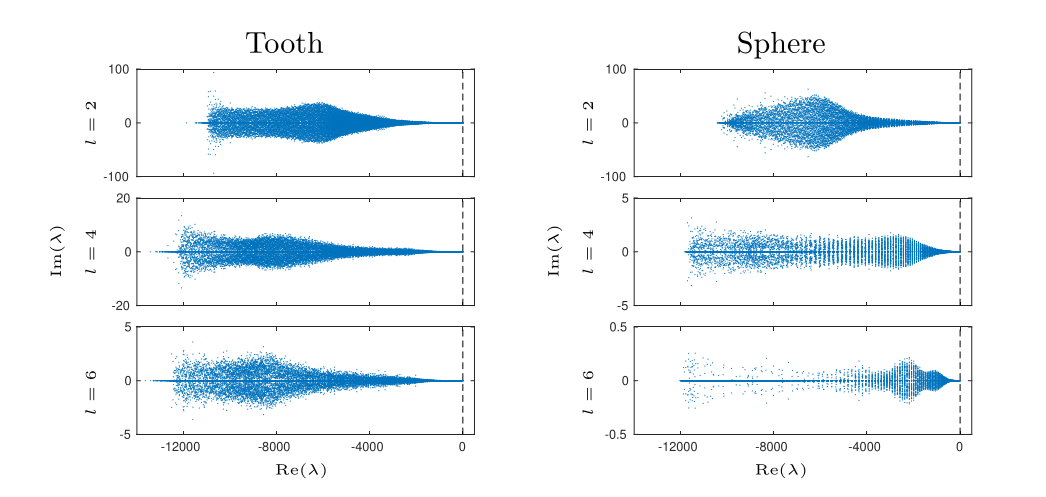}
\caption{Eigenvalue spectra of $\Delta_{\Gamma}$ for the tooth and unit sphere from Section \ref{sec_param3Dcase} with $N = 20{,}298$ ($h = 0.036$) and $N = 10{,}000$ ($h = 0.035$) nodes, respectively, when using RBF-FD $r^5$ augmented with polynomials of degree $l=2,4$ and $6$.}\label{Fig:eigen_ToothSphere}
\end{figure}

Figure \ref{Fig:eigen_ToothSphere} shows the eigenvalue spectra of the differentiation matrices from Section \ref{sec_param3Dcase} approximating $\Delta_{\Gamma}$ on the tooth and unit sphere with $N = 20{,}298$ ($h = 0.036$) and $N = 10{,}000$ ($h = 0.035$) nodes, when using RBF-FD $r^5$ augmented with polynomials of degree $l=2,4$ and 6. As the polynomial degree increases, the differentiation matrix provides a better approximation of the surface Laplacian, with the imaginary part of the eigenvalues approaching zero (see the zoom on the $y$-axis). Note that none of the eigenvalues fall in the right half-plane.

\subsection{Further results on Turing patterns} \label{App:Teig}

\begin{figure}
\centering

\includegraphics[scale=0.35,trim={2ex 0ex 0 3ex},clip]{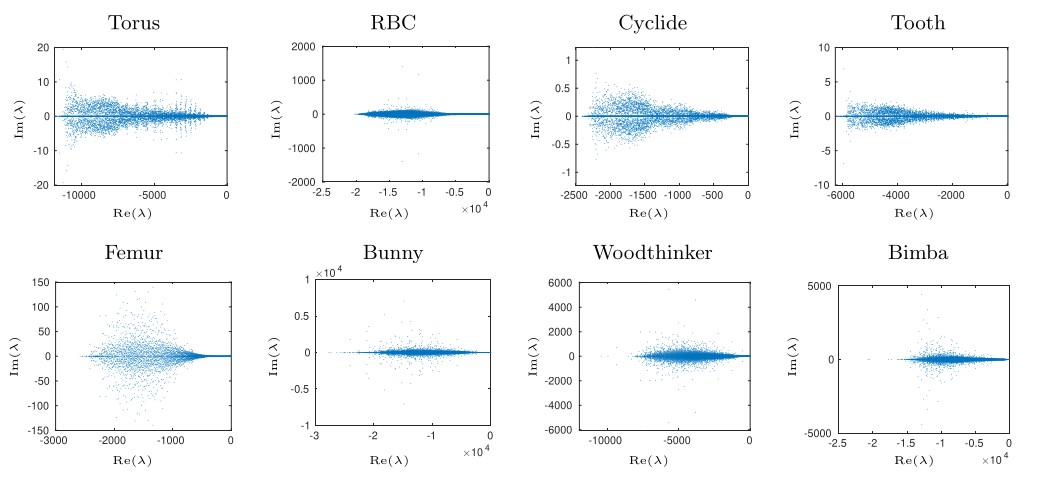}

\caption{Eigenvalues of differentiation matrix approximating the Laplace--Beltrami operator on the surface models from Figure \ref{Fig:Turing}. Observe that none of them fall in the right halfplane.}
\label{Fig:eigen}
\end{figure}

Figure~\ref{Fig:eigen} shows the eigenvalue distributions of the Laplace--Beltrami operator on the surfaces considered in Section~\ref{sec_TuringStatic}. Recall that PHS $r^5$ combined with polynomials of degree $6$ are employed. In all the cases, stable approximations with negative real part eigenvalues are obtained.

\section{Code availability}
Sample MATLAB codes that reproduce some of the results presented in this paper are available upon request at \texttt{\url{https://gitlab.com/apetras/surfacepdes_public_code.git}}.

\newpage
\bibliographystyle{siamplain}
\bibliography{draftBibliography}

\begin{thebibliography}{10}

\bibitem{ascher1995implicit}
{\sc U.~M. Ascher, S.~J. Ruuth, and B.~T. Wetton}, {\em Implicit-explicit
  methods for time-dependent partial differential equations}, SIAM Journal on
  Numerical Analysis, 32 (1995), pp.~797--823,
  \url{https://doi.org/10.1137/0732037}.

\bibitem{B19}
{\sc V.~Bayona}, {\em An insight into {RBF-FD} approximations augmented with
  polynomials}, Computers \& Mathematics with Applications, 77 (2019),
  pp.~2337--2353, \url{https://doi.org/10.1016/j.camwa.2018.12.029}.

\bibitem{BFFB17}
{\sc V.~Bayona, N.~Flyer, B.~Fornberg, and G.~A. Barnett}, {\em On the role of
  polynomials in {RBF-FD} approximations: {II. N}umerical solution of elliptic
  {PDE}s}, J. Comput. Phys., 332 (2017), pp.~257--273,
  \url{https://doi.org/10.1016/j.jcp.2016.12.008}.

\bibitem{chen2019meshless}
{\sc M.~Chen, K.~C. Cheung, and L.~LING}, {\em {Meshless collocation methods
  for solving PDES on surfaces}}, WIT Transactions on Engineering Sciences, 126
  (2019), pp.~159--170, \url{https://doi.org/10.2495/BE420141}.

\bibitem{chen2023kernel}
{\sc M.~Chen, K.~C. Cheung, and L.~Ling}, {\em A kernel-based least-squares
  collocation method for surface diffusion}, SIAM Journal on Numerical
  Analysis, 61 (2023), pp.~1386--1404,
  \url{https://doi.org/10.1137/21M1444369}.

\bibitem{chen2020extrinsic}
{\sc M.~Chen and L.~Ling}, {\em Extrinsic meshless collocation methods for
  {PDEs} on manifolds}, SIAM Journal on Numerical Analysis, 58 (2020),
  pp.~988--1007, \url{https://doi.org/10.1137/17M1158641}.

\bibitem{chen2020kernel}
{\sc M.~Chen and L.~Ling}, {\em Kernel-based collocation methods for heat
  transport on evolving surfaces}, Journal of Computational Physics, 405
  (2020), p.~109166, \url{https://doi.org/10.1016/j.jcp.2019.109166}.

\bibitem{cheung2018kernel}
{\sc K.~C. Cheung and L.~Ling}, {\em {A kernel-based embedding method and
  convergence analysis for surfaces PDEs}}, SIAM Journal on Scientific
  Computing, 40 (2018), pp.~A266--A287,
  \url{https://doi.org/10.1137/16M1080410}.

\bibitem{cheung2015localized}
{\sc K.~C. Cheung, L.~Ling, and S.~J. Ruuth}, {\em A localized meshless method
  for diffusion on folded surfaces}, Journal of Computational Physics, 297
  (2015), pp.~194--206, \url{https://doi.org/10.1016/j.jcp.2015.05.021}.

\bibitem{demanet2006painless}
{\sc L.~Demanet}, {\em {Painless, highly accurate discretizations of the
  Laplacian on a smooth manifold}}, tech. report, Citeseer, 2006.

\bibitem{Dziuk2013}
{\sc G.~Dziuk and C.~M. Elliott}, {\em {Finite element methods for surface
  PDEs}}, Acta Numerica, 22 (2013), pp.~289--396,
  \url{https://doi.org/10.1017/S0962492913000056}.

\bibitem{engwirda2016off}
{\sc D.~Engwirda and D.~Ivers}, {\em Off-centre {S}teiner points for
  {D}elaunay-refinement on curved surfaces}, Computer-Aided Design, 72 (2016),
  pp.~157--171, \url{https://doi.org/10.1016/j.cad.2015.10.007}.

\bibitem{fasshauer2007meshfree}
{\sc G.~E. Fasshauer}, {\em {Meshfree approximation methods with MATLAB}},
  vol.~6, World Scientific, 2007, \url{https://doi.org/10.1142/6437}.

\bibitem{FFBB16}
{\sc N.~Flyer, B.~Fornberg, V.~Bayona, and G.~A. Barnett}, {\em On the role of
  polynomials in {RBF-FD} approximations: {I. I}nterpolation and accuracy}, J.
  Comput. Phys., 321 (2016), pp.~21--38,
  \url{https://doi.org/10.1016/j.jcp.2016.05.026}.

\bibitem{flyer2012guide}
{\sc N.~Flyer, E.~Lehto, S.~Blaise, G.~B. Wright, and A.~St-Cyr}, {\em {A guide
  to RBF-generated finite differences for nonlinear transport: Shallow water
  simulations on a sphere}}, Journal of Computational Physics, 231 (2012),
  pp.~4078--4095, \url{https://doi.org/10.1016/j.jcp.2012.01.028}.

\bibitem{flyer2009radial}
{\sc N.~Flyer and G.~B. Wright}, {\em A radial basis function method for the
  shallow water equations on a sphere}, Proceedings of the Royal Society A:
  Mathematical, Physical and Engineering Sciences, 465 (2009), pp.~1949--1976,
  \url{https://doi.org/10.1098/rspa.2009.0033}.

\bibitem{fornberg2015primer}
{\sc B.~Fornberg and N.~Flyer}, {\em A primer on radial basis functions with
  applications to the geosciences}, SIAM, 2015.

\bibitem{fornberg_flyer_2015}
{\sc B.~Fornberg and N.~Flyer}, {\em {Solving PDEs with radial basis
  functions}}, Acta Numerica, 24 (2015), pp.~215--258,
  \url{https://doi.org/10.1017/S0962492914000130}.

\bibitem{fornberg2011stabilization}
{\sc B.~Fornberg and E.~Lehto}, {\em Stabilization of {RBF}-generated finite
  difference methods for convective {PDE}s}, Journal of Computational Physics,
  230 (2011), pp.~2270--2285.

\bibitem{fuselier2013high}
{\sc E.~J. Fuselier and G.~B. Wright}, {\em A high-order kernel method for
  diffusion and reaction-diffusion equations on surfaces}, Journal of
  Scientific Computing, 56 (2013), pp.~535--565,
  \url{https://doi.org/10.1007/s10915-013-9688-x}.

\bibitem{gunderman2020transport}
{\sc D.~Gunderman, N.~Flyer, and B.~Fornberg}, {\em {Transport schemes in
  spherical geometries using spline-based RBF-FD with polynomials}}, Journal of
  Computational Physics, 408 (2020), p.~109256,
  \url{https://doi.org/10.1016/j.jcp.2020.109256}.

\bibitem{jones2023generalized}
{\sc A.~M. Jones, P.~A. Bosler, P.~A. Kuberry, and G.~B. Wright}, {\em
  Generalized moving least squares vs. radial basis function finite difference
  methods for approximating surface derivatives}, Computers \& Mathematics with
  Applications, 147 (2023), pp.~1--13,
  \url{https://doi.org/10.1016/j.camwa.2023.07.015}.

\bibitem{lehto2017radial}
{\sc E.~Lehto, V.~Shankar, and G.~B. Wright}, {\em A radial basis function
  {(RBF)} compact finite difference {(FD)} scheme for reaction-diffusion
  equations on surfaces}, SIAM Journal on Scientific Computing, 39 (2017),
  pp.~A2129--A2151, \url{https://doi.org/10.1137/16M1095457}.

\bibitem{leung2009grid}
{\sc S.~Leung and H.~Zhao}, {\em A grid based particle method for moving
  interface problems}, Journal of Computational Physics, 228 (2009),
  pp.~2993--3024, \url{https://doi.org/10.1016/j.jcp.2009.01.005}.

\bibitem{macdonald2009implicit}
{\sc C.~B. Macdonald and S.~J. Ruuth}, {\em The implicit closest point method
  for the numerical solution of partial differential equations on surfaces},
  SIAM Journal on Scientific Computing, 31 (2009), pp.~4330--4350,
  \url{https://doi.org/10.1137/080740003}.

\bibitem{marz2012calculus}
{\sc T.~Marz and C.~B. Macdonald}, {\em Calculus on surfaces with general
  closest point functions}, SIAM Journal on Numerical Analysis, 50 (2012),
  pp.~3303--3328, \url{https://doi.org/10.1137/120865537}.

\bibitem{petras2019least}
{\sc A.~Petras, L.~Ling, C.~Piret, and S.~J. Ruuth}, {\em A least-squares
  implicit {RBF-FD} closest point method and applications to {PDEs} on moving
  surfaces}, Journal of Computational Physics, 381 (2019), pp.~146--161,
  \url{https://doi.org/10.1016/j.jcp.2018.12.031}.

\bibitem{petras2018rbf}
{\sc A.~Petras, L.~Ling, and S.~J. Ruuth}, {\em An {RBF-FD} closest point
  method for solving {PDEs} on surfaces}, Journal of Computational Physics, 370
  (2018), pp.~43--57, \url{https://doi.org/10.1016/j.jcp.2018.05.022}.

\bibitem{petras2022meshfree}
{\sc A.~Petras, L.~Ling, and S.~J. Ruuth}, {\em {Meshfree Semi-Lagrangian
  methods for solving surface advection PDEs}}, Journal of Scientific
  Computing, 93 (2022), pp.~1--22,
  \url{https://doi.org/10.1007/s10915-022-01966-w}.

\bibitem{petras2016pdes}
{\sc A.~Petras and S.~J. Ruuth}, {\em {PDEs on moving surfaces via the closest
  point method and a modified grid based particle method}}, Journal of
  Computational Physics, 312 (2016), pp.~139--156,
  \url{https://doi.org/10.1016/j.jcp.2016.02.024}.

\bibitem{piret2012orthogonal}
{\sc C.~Piret}, {\em {The orthogonal gradients method: A radial basis functions
  method for solving partial differential equations on arbitrary surfaces}},
  Journal of Computational Physics, 231 (2012), pp.~4662--4675,
  \url{https://doi.org/10.1016/j.jcp.2012.03.007}.

\bibitem{piret2016fast}
{\sc C.~Piret and J.~Dunn}, {\em {Fast RBF OGr for solving PDEs on arbitrary
  surfaces}}, in AIP Conference Proceedings, vol.~1776, AIP Publishing, 2016,
  \url{https://doi.org/10.1063/1.4965351}.

\bibitem{reeger2016numerical}
{\sc J.~A. Reeger and B.~Fornberg}, {\em Numerical quadrature over the surface
  of a sphere}, Studies in Applied Mathematics, 137 (2016), pp.~174--188,
  \url{https://doi.org/10.1111/sapm.12106}.

\bibitem{reeger2018numerical}
{\sc J.~A. Reeger and B.~Fornberg}, {\em Numerical quadrature over smooth
  surfaces with boundaries}, Journal of Computational Physics, 355 (2018),
  pp.~176--190, \url{https://doi.org/10.1016/j.jcp.2017.11.010}.

\bibitem{ruuth2008simple}
{\sc S.~J. Ruuth and B.~Merriman}, {\em A simple embedding method for solving
  partial differential equations on surfaces}, Journal of Computational
  Physics, 227 (2008), pp.~1943--1961,
  \url{https://doi.org/10.1016/j.jcp.2007.10.009}.

\bibitem{shankar2018rbf}
{\sc V.~Shankar, A.~Narayan, and R.~M. Kirby}, {\em {RBF-LOI: Augmenting radial
  basis functions (RBFs) with least orthogonal interpolation (LOI) for solving
  PDEs on surfaces}}, Journal of Computational Physics, 373 (2018),
  pp.~722--735, \url{https://doi.org/10.1016/j.jcp.2018.07.015}.

\bibitem{shankar2015radial}
{\sc V.~Shankar, G.~B. Wright, R.~M. Kirby, and A.~L. Fogelson}, {\em A radial
  basis function {(RBF)}-finite difference {(FD)} method for diffusion and
  reaction--diffusion equations on surfaces}, Journal of Scientific Computing,
  63 (2015), pp.~745--768, \url{https://doi.org/10.1007/s10915-014-9914-1}.

\bibitem{shankar2020robust}
{\sc V.~Shankar, G.~B. Wright, and A.~Narayan}, {\em {A robust hyperviscosity
  formulation for stable RBF-FD discretizations of advection-diffusion-reaction
  equations on manifolds}}, SIAM Journal on Scientific Computing, 42 (2020),
  pp.~A2371--A2401, \url{https://doi.org/10.1137/19M1288747}.

\bibitem{shaw2019radial}
{\sc S.~B. Shaw}, {\em {Radial basis function finite difference approximations
  of the Laplace-Beltrami operator}},  (2019),
  \url{https://doi.org/10.18122/td/1587/boisestate}.

\bibitem{suchde2019fully}
{\sc P.~Suchde and J.~Kuhnert}, {\em {A fully Lagrangian meshfree framework for
  PDEs on evolving surfaces}}, Journal of Computational Physics, 395 (2019),
  pp.~38--59, \url{https://doi.org/10.1016/j.jcp.2019.06.031}.

\bibitem{Suchde2019}
{\sc P.~Suchde and J.~Kuhnert}, {\em {A meshfree generalized finite difference
  method for surface PDEs}}, Computers {\&} Mathematics with Applications, 78
  (2019), pp.~2789--2805, \url{https://doi.org/10.1016/j.camwa.2019.04.030},
  \url{https://arxiv.org/abs/1806.07193}.

\bibitem{tang2021localized}
{\sc Z.~Tang, Z.~Fu, M.~Chen, and L.~Ling}, {\em {A localized extrinsic
  collocation method for Turing pattern formations on surfaces}}, Applied
  Mathematics Letters, 122 (2021), p.~107534,
  \url{https://doi.org/10.1016/j.aml.2021.107534}.

\bibitem{von2013embedded}
{\sc I.~von Glehn, T.~M{\"a}rz, and C.~B. Macdonald}, {\em An embedded
  method-of-lines approach to solving partial differential equations on
  surfaces}, arXiv preprint arXiv:1307.5657,  (2013),
  \url{https://doi.org/10.48550/arXiv.1307.5657}.

\end{thebibliography}
\end{document}